%% file: revise1.tex
\documentclass[final]{siamltex}		
\usepackage{verbatim,algorithmic} 

\usepackage{color,url,hyperref}

\usepackage{graphicx,multirow}  
\usepackage{amssymb,amsmath}

\input{preamble.tex}

\usepackage{graphicx,amsmath,amsfonts,amssymb,subfig}
\usepackage[ruled,vlined,linesnumbered]{algorithm2e}

\definecolor{dred}{rgb}{0.6,0.0,0.0}

\newcommand{\tphi}{\wt{\phi}}
\newcommand{\gsz}{g_{\sigma}}

\usepackage{placeins}
\graphicspath{{FIGS072014/}}

\title{Approximating spectral densities of large matrices}

\author{Lin Lin 
\thanks{Department of Mathematics, University of California, Berkeley
and Computational Research Division, Lawrence Berkeley National
Laboratory, Berkeley, CA 94720 \texttt{linlin@math.berkeley.edu}} 
\and Yousef Saad 
\thanks{Department of Computer Science and Engineering,
University of Minnesota, Twin Cities, Minneapolis, MN 55455 \texttt{saad@cs.umn.edu}}
\and Chao Yang \thanks{Computational Research Division, Lawrence Berkeley National Laboratory, Berkeley, CA 94720 \texttt{cyang@lbl.gov}} 
}

\begin{document} 

\maketitle 

\begin{abstract}
  In physics, it is sometimes desirable to compute
  the  so-called \emph{Density  Of States}  (DOS), also  known  as the
  \emph{spectral  density},  of  a  real symmetric  matrix
  $A$. The  spectral density  can be viewed  as a  probability density
  distribution  that measures  the likelihood  of  finding eigenvalues
  near some point  on the real line.  The  most straightforward way to
  obtain this density  is to compute all eigenvalues  of $A$. But this
  approach is  generally costly and wasteful,  especially for matrices
  of large  dimension.  There exists alternative methods  that allow us to
  estimate the  spectral density function at  much lower cost.  The
  major computational cost of these methods is in multiplying $A$ with
  a number  of vectors, which makes them appealing for
  large-scale  problems  where products  of  the  matrix $A$  with
  arbitrary vectors  are relatively inexpensive.  This paper defines
  the problem of estimating the spectral  density carefully, and
  discusses how to measure the accuracy of an approximate spectral
  density.  It then surveys a  few  known  methods  for  estimating the
  spectral density,  and proposes some  new variations of existing
  methods.  All  methods are  discussed from a  numerical linear algebra
  point of view.
\end{abstract} 


 \pagestyle{myheadings}
\thispagestyle{plain}
\markboth{L. LIN, Y. SAAD AND C. YANG}{SPECTRAL DENSITIES OF LARGE MATRICES}

\section{Introduction} \label{sec:intro}
Given an $n\times n$ real symmetric and sparse matrix $A$, scientists in various
disciplines often want to compute its \emph{Density Of States} (DOS), 
or \emph{spectral density}. 
Formally, the DOS is defined as
\begin{equation}
  \phi(t) = \frac{1}{n}  \sum_{j=1}^n \delta(t - \lambda_j), 
  \label{eq:DOS0}
\end{equation}
where $\delta$ is the Dirac distribution commonly referred to as
the  Dirac $\delta$-``function''~\cite{Schwartz1966,ByronFuller,Richtmyer1981}, 
and the
$\lambda_j$'s are the eigenvalues of $A$, assumed here to be labeled
non-decreasingly.  
Using the DOS, the number of eigenvalues in an
interval $[a,b]$ can be formally expressed as
\begin{equation}
  \label{eq:DOS1}
  \nu_{[a, b]} = \int_a^b 
  \sum_j  \delta(t - \lambda_j) \ dt \equiv 
  \int_a^b n \phi(t) dt.
\end{equation}
Therefore, one can view $\phi(t)$ as a probability
distribution ``function'', which gives the probability of finding 
eigenvalues of $A$ in a given infinitesimal interval near $t$. 
If one had access to all the eigenvalues of $A$, the task of 
computing the DOS would become a trivial one.  However, in many  applications, 
the dimension of $A$ is large. The computation of its entire spectrum
is prohibitively expensive, and this leads to the need to develop 
efficient alternative methods to
estimate $\phi(t)$ without computing eigenvalues of $A$.
Since $\phi(t)$ is not a proper function, we need to  clarify what we mean
by ``estimating" $\phi(t)$, and this will be addressed in detail shortly.
For now we can use our intuition to argue  that
 $\phi(t)$ can be approximated by dividing the interval  containing the 
spectrum  of $A$ into many sub-intervals and use a tool like Sylvester's
law of inertia to count the number of eigenvalues within each of these
sub-intervals.
This approach yields a histogram of the eigenvalues. 
Expression \eqref{eq:DOS1} will then provide us with an average
``value'' of 
$\phi(t)$ in each small subinterval $[a, b]$. 
As the size of
each subinterval decreases, the histogram approaches the spectral
density of $A$. However,  this is not a practical approach
since performing such an inertia count requires us 
to compute the $LDL^T$ factorization \cite{GVL-book} of 
$A - t_iI$, where the $t_i$'s are the end points of the subintervals.
In general, this approach is prohibitively expensive because of the large
 number of intervals needed and the shear cost of each factorization.
Therefore, a procedure that relies entirely on multiplications
 of $A$ with vectors is the only viable  approach. 

Because calculating the spectral density is such an important problem in
quantum mechanics, there is an abundant literature devoted to this
problem and research in this area was extremely active in the 1970s and
1980s, leading to  clever and powerful methods  developed by physicists
and chemists~\cite{Ducastelle1970,Turek88,drobold93,wheeler72} for this
purpose.

In this survey paper, we review two classes of methods
for approximating the spectral density of a real symmetric matrix 
from a numerical linear algebra perspective.  
For simplicity, all methods are presented using real arithmetic
operations, i.e. we assume the matrix is real symmetric. The
generalization to Hermitian matrices is straightforward.
The first class of methods contains the Kernel Polynomial
Method (KPM)~\cite{SilverRoder1994,Wang1994} and its variants.  
The KPM can be viewed as a formal polynomial expansion of the 
spectral density.  It uses a moment matching method to derive 
the coefficients for the polynomials. 
The method, which is widely used in a variety of 
calculations that require the DOS \cite{kpmsurvey2006},  has continued
to receive a tremendous amount of interest in the last few years 
\cite{kpmsurvey2006,cpb2010,jck2012,spd2013}. 
We show that a less well known, but rather original, method due to
Lanczos and known as the ``Lanczos spectroscopic'' procedure, 
which samples the cosine transform of the spectral density, 
is closely related to KPM. As another variant of KPM, we 
present a spectral density probing method called Delta-Gauss-Legendre method,
that can be viewed as a polynomial expansion of a smoothed 
spectral density.
The second class of methods we consider uses the classical 
Lanczos procedure to partially diagonalize $A$. The 
eigenvalues and eigenvectors of tridiagonal matrix are
used to construct approximations to the spectral density.

One of the key ingredients used in most of these methods is a 
well-established artifice for estimating the trace of a matrix. For example,
the expansion coefficients in the above-mentioned KPM method can be obtained
from the traces of the matrix polynomials $T_k(A)$, where $T_k$ is 
the Chebyshev polynomial of degree $k$.
 Each of these is in turn
estimated as the mean of $v^T T_k(A) v$ over a number of random vectors $v$.
This procedure for estimating the trace has been discovered more or less
independently by statisticians \cite{Hutchinson1989} and
physicists and chemists \cite{SilverRoder1994,Wang1994}. 

A natural question one would ask is: among all the methods reviewed
here, which is the best method to use? The answer to this question is not 
simple.  
Since the methods discussed in this paper are all based on 
matrix-vector product operations (MATVECs) 
the criterion for choosing the best method should be based on the 
quality of the approximation when approximately the same
number of MATVECs are used. 
In order to determine 
the quality, we must first establish a way to measure 
the accuracy of the approximation.
Because the spectral density is defined in terms of the 
Dirac $\delta$-``function"s, which are not proper functions but 
distributions~\cite{ByronFuller,Richtmyer1981}, 
the standard error metrics used for approximating smooth functions 
are not appropriate.  
Furthermore, the accuracy measure 
should  depend on the desired resolution. In many applications,
it is not necessary to obtain a high resolution spectral density.
If fact such a high resolution density would be  highly
discontinuous, if one thinks of our already mentioned intuitive
interpretation in terms of a histogram.
For these reasons our proposed  metric for measuring the accuracy of spectral
density approximation is  defined in section~\ref{sec:dosFun}, so as to
 allow rigorous quantitative comparisons of the 
different spectral density approximation methods, 
instead of  relying on a subjective visual measure as is often
done in practice.

All approximation methods we consider are presented in section~\ref{sec:methods}. 
We give some numerical examples in section~\ref{sec:examples} 
to compare different numerical methods for approximating
spectral densities.  We illustrate the effectiveness of our error
metric for evaluating the quality of the approximation, and 
describe some general observations we made about the behavior
of different methods.

\section{Assessing the quality of spectral density approximation}\label{sec:dosFun}
We will denote the approximate spectral density by $\tphi(t)$ which is 
a regular function.
The types of approximate spectral densities we consider in this paper are
all continuous functions. However, 
since $\phi(t)$ is defined in terms of a number of Dirac $\delta$-functions
that are not proper functions but distributions, we cannot use 
the standard $L^p$-norm, with, e.g. $p=1,2$, or $\infty$, to evaluate
the approximation error defined in terms of $\phi(t)-\tphi(t)$,
where we note
 that in this  difference $\tphi$ is interpreted as a distribution.

We discuss two approaches to get around this difficulty.
In the first approach, we use the fact that $\delta(t)$ is a
distribution,
i.e. it is formally defined 
through applications to a test function $g$: 
\[
\langle\delta(\cdot-\lambda), g\rangle \equiv \int_{-\infty}^{\infty}
\delta(t-\lambda) g (t) dt \equiv g(\lambda), 
\] 
where we use $\delta(\cdot-\lambda)$ to denote a Dirac $\delta$ centered 
at $\lambda$, $g\in C^{\infty}(\mathbb{R})$, and for all $p,k\in
\mathbb{N}$,
\[
\sup_{t\in \mathbb{R}} \abs{t^{p} g^{(k)}(t)} < \infty.
\]
Here $g^{(k)}(t)$ is the $k$th derivative of $g(t)$.  The test function
$g$ is chosen to be a member of the Schwartz space (or Schwartz
class)~\cite{Richtmyer1981}, denoted by $\mathcal{S}$. In other words,
the test function $g$ should be smooth and decays sufficiently fast 
towards 0 when $|t|$ approaches infinity.
The error is then measured as
\begin{equation}
\epsilon_1 = \sup_{g\in \mathcal{S}} \abs{\average{\phi,g} - \average{\wt{\phi},g}}.
\label{eq:err1}
\end{equation}
In practice, we restrict $\mathcal{S}$ to be a subspace of the Schwartz space
that allows us to compute~\eqref{eq:err1} at a finite resolution.
We will elaborate on the choice of $g$ and $\mathcal{S}$ in section~\ref{sec:err1}.

In the second approach, we regularize $\delta$-functions and replace them
with continuous and smooth functions such as Gaussians with an appropriately
chosen standard deviation $\sigma$.  The resulting regularized
spectral density, which we denote by $\phi_{\sigma}(t)$, is a well defined
function. Hence, it is meaningful to compute the approximation error
\begin{equation}
\epsilon_2 = \|\phi_{\sigma} (t)-\tphi(t)\|_p,
\label{eq:err2}
\end{equation}
for $p=1,2$, and $\infty$. 
There is a close connection between the first and second approach, on which 
we will elaborate in the next section.

We should note that the notion of regularization, which is rarely 
discussed in the existing physics and chemistry literature, is important for assessing 
the accuracy of spectral density approximation.  A fully accurate
approximation amounts to computing all eigenvalues of $A$. But
for most applications, one only needs to know the number of 
eigenvalues within any small subinterval contained in the spectrum of $A$. 
The size of the interval represents the ``resolution" of the approximation.  
The accuracy of the approximation is only meaningful up to the 
desired resolution.  When~\eqref{eq:err2} is used to assess the quality 
of the approximation, the resolution is defined in terms of the regularization 
parameter $\sigma$.  A smaller $\sigma$ corresponds to higher resolution.

The notion of resolution can also be built into the error metric
~\eqref{eq:err1} if the trial function $g$ belongs to 
a certain class of functions, which we will discuss in 
the next section.

\subsection{Restricting the test function space $\mathcal{S}$ }\label{sec:err1}

The fact that the spectral density $\phi(t)$ is defined in terms of 
Dirac $\delta$-functions suggests that no smooth function can 
converge to the spectral density in the limit of high resolution.

To see this, consider  $\nu_{[a,b]}$ defined in Eq.~\eqref{eq:DOS1} and the
associated approximation obtained from a smooth approximation $\wt{\phi}(t)$ as
\[
\wt{\nu}_{[a,b]} = \int_{a}^{b} n \wt{\phi}(t) dt.
\]
For simplicity, let the spectral density $\phi(t)=\delta(t)$ to be a
single $\delta$-function, and the number of eigenvalues $n=1$. Infinite
resolution means that $\abs{\nu_{[a,b]}-\wt{\nu}_{[a,b]}}$ should be
small for any choice of $[a,b]$.  Now take $a=-\varepsilon,b=\varepsilon$. It is
easy to verify that 
\[
\lim_{\varepsilon\to 0+}\nu_{[-\varepsilon,\varepsilon]} = 1, \quad
\lim_{\varepsilon\to 0+}\wt{\nu}_{[-\varepsilon,\varepsilon]} = 0,
\]
In this sense, \textit{all} smooth approximation of the spectral density
results in the same accuracy, i.e. there is no difference between a
carefully designed approximation of the spectral density and a constant
approximation.  
Hence, the distribution $\phi(t)$ behaves very much
like a highly discontinuous function and cannot be approximated by
smooth functions with infinite resolution.


In practice, physical quantities and observables can often be deduced 
from spectral density at finite resolution, i.e. the eigenvalue count 
only needs to be approximately correct for an interval of a given finite size.  
For instance, in condensed matter physics, such
information is enough to provide material properties such as the band
gap or the Van Hove singularity~\cite{AshcroftMermin1976} within given target accuracy.
The reduced resolution requirement suggests that we may not need to 
take the test space $S$ in~\eqref{eq:err1} to be the whole Schwartz space.
Instead, we can choose functions that have ``limited resolution" as 
test functions.  For example, we may consider using Gaussian functions 
of the form
\begin{equation}
  g_{\sigma}(t) = \frac{1}{(2\pi \sigma^2)^{1/2}} e^{-\frac{t^2}{2
  \sigma^2} },
  \label{eq:gauss}
\end{equation}
and restrict $\mathcal{S}$ to the subspace 
\[
\mathcal{S}(\sigma;[\lambda_{lb},\lambda_{ub}]) = \left\{g \Big\vert g(t) \equiv
g_{\sigma}(t-\lambda), \quad \lambda\in [\lambda_{lb},\lambda_{ub}]
\right\},
\]
where $\lambda_{lb}$ and $\lambda_{ub}$ are lower and upper bounds of the eigenvalues of 
$A$ respectively, and the parameter $\sigma$ defines the
\textit{target resolution} up to which we intend to measure. 
The use of Gaussian functions in the space $\mathcal{S}(\sigma;[\lambda_{lb},\lambda_{ub}])$ 
can be understood as a smooth way of counting the number of eigenvalues 
in an interval whose size is proportional to $\sigma$.

Using this choice of the test space, we can measure the quality of any approximation 
by the following metric:
\begin{equation}
  E[\tilde{\phi}; \mathcal{S}(\sigma;[\lambda_{lb},\lambda_{ub}])] = \sup_{g\in
  \mathcal{S}(\sigma;[\lambda_{lb},\lambda_{ub}])} 
  \abs{\average{\phi, g} - 
  \average{\tilde{\phi}, g}}.
  \label{eqn:errormetric}
\end{equation}

We remark that the use of Gaussians is not the only way to
restrict the test space. In some applications, the DOS is often used 
as a measure for integrating certain physical quantities of interest. 
If the quantity of interest can be expressed as
\[
\average{\phi,g} \equiv \int g(\lambda) \phi(\lambda) d\lambda \equiv 
\frac{1}{n}\sum_{j=1}^{n} g(\lambda_{j})
\]
for some smooth function $g$, then $\mathcal{S}$ can be chosen to
contain only one function $g$, and the approximation error is naturally defined as
\begin{equation}
  E[\tilde{\phi}; g] =  \abs{\average{\phi(t), g} - 
  \average{\tilde{\phi}(t), g}},
  \label{eqn:errormetric_g}
\end{equation}
for that particular function $g$. We give an example of
this measure in section~\ref{sec:application}.

\subsection{Regularizing the spectral density} \label{sec:reg}
The error metric in Eq.~\eqref{eqn:errormetric} can also
be understood in the following sense. Let
\begin{equation}
\phi_{\sigma}(t) = \average{\phi(\cdot), g_{\sigma}(\cdot-t)} =
\sum_{j=1}^{n} g_{\sigma}(t-\lambda_{j}),
\label{eq:blurphi}
\end{equation}
then $\phi_{\sigma}(t)$ is nothing but a blurred or regularized spectral
density, and the blurring is given by a Gaussian function with width $\sigma$.  Similarly,
$\average{\wt{\phi}(\cdot), \gsz(\cdot-t)}$
can be understood as a blurred version of an approximate spectral density.
Therefore the error metric in Eq.~\eqref{eqn:errormetric} is equivalent
to the $L^{\infty}$ error between two well defined functions
$\average{\wt{\phi}(\cdot), \gsz(\cdot-t)}$ against $\phi_{\sigma}(t)$. 

This point of view leads to another way to construct and 
measure the approximation to the spectral density function.
Instead of trying to approximate $\phi(t)$ directly, which 
may be difficult due to the presence of $\delta$-functions in $\phi(t)$,
we first construct a smooth representation of the $\delta$-function.
The representation we choose should be commensurate with the 
desired resolution of the spectral density.
This regularization process allows us to expand smooth functions
in terms of other smooth functions such as orthogonal polynomials, and
the approximation error associated with such expansions can
be evaluated directly and without introducing additional
regularization procedure.

\begin{figure}[htb] 
\begin{center}
\includegraphics[width=0.30\textwidth]{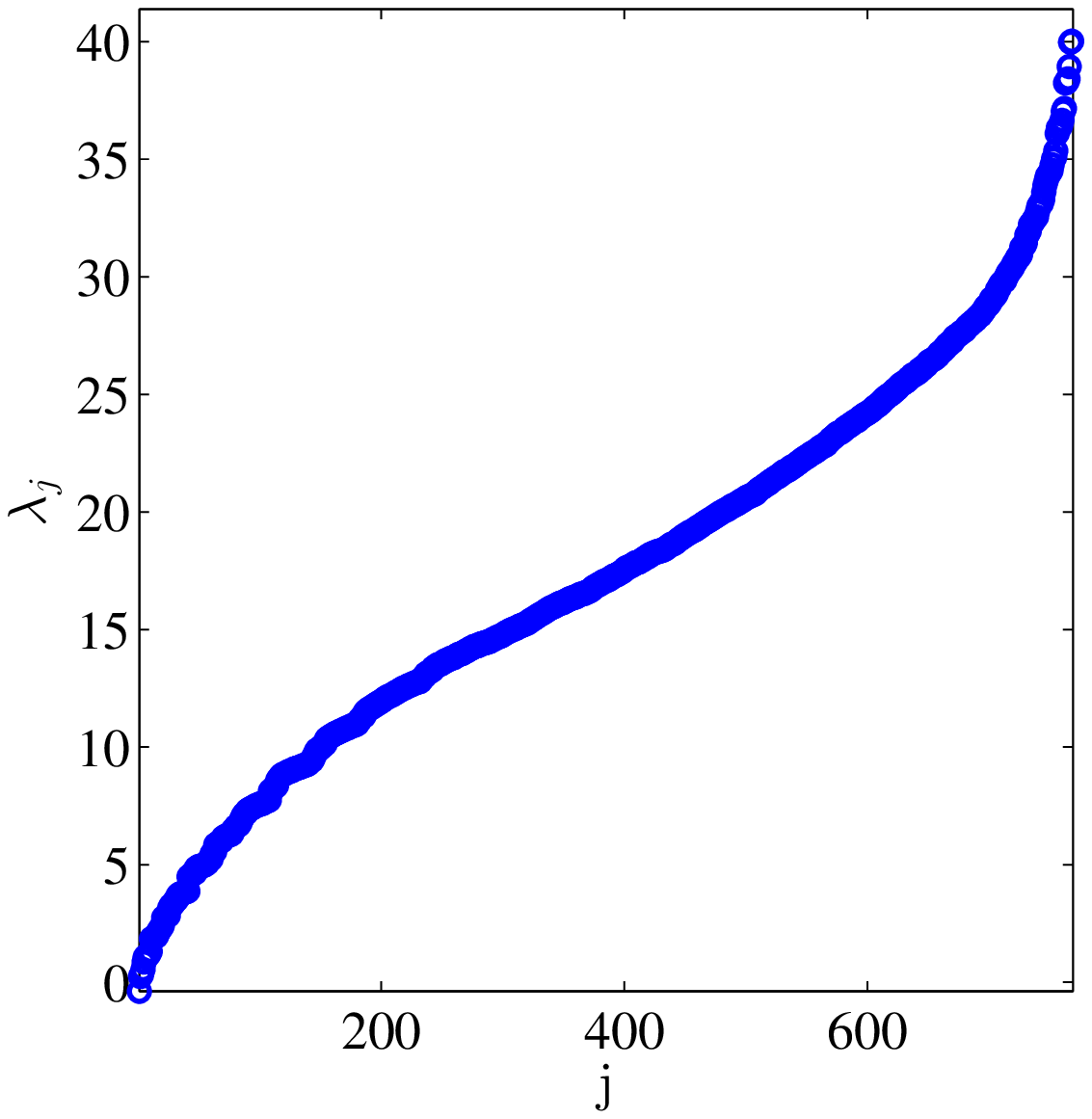}
\quad
\includegraphics[width=0.30\textwidth]{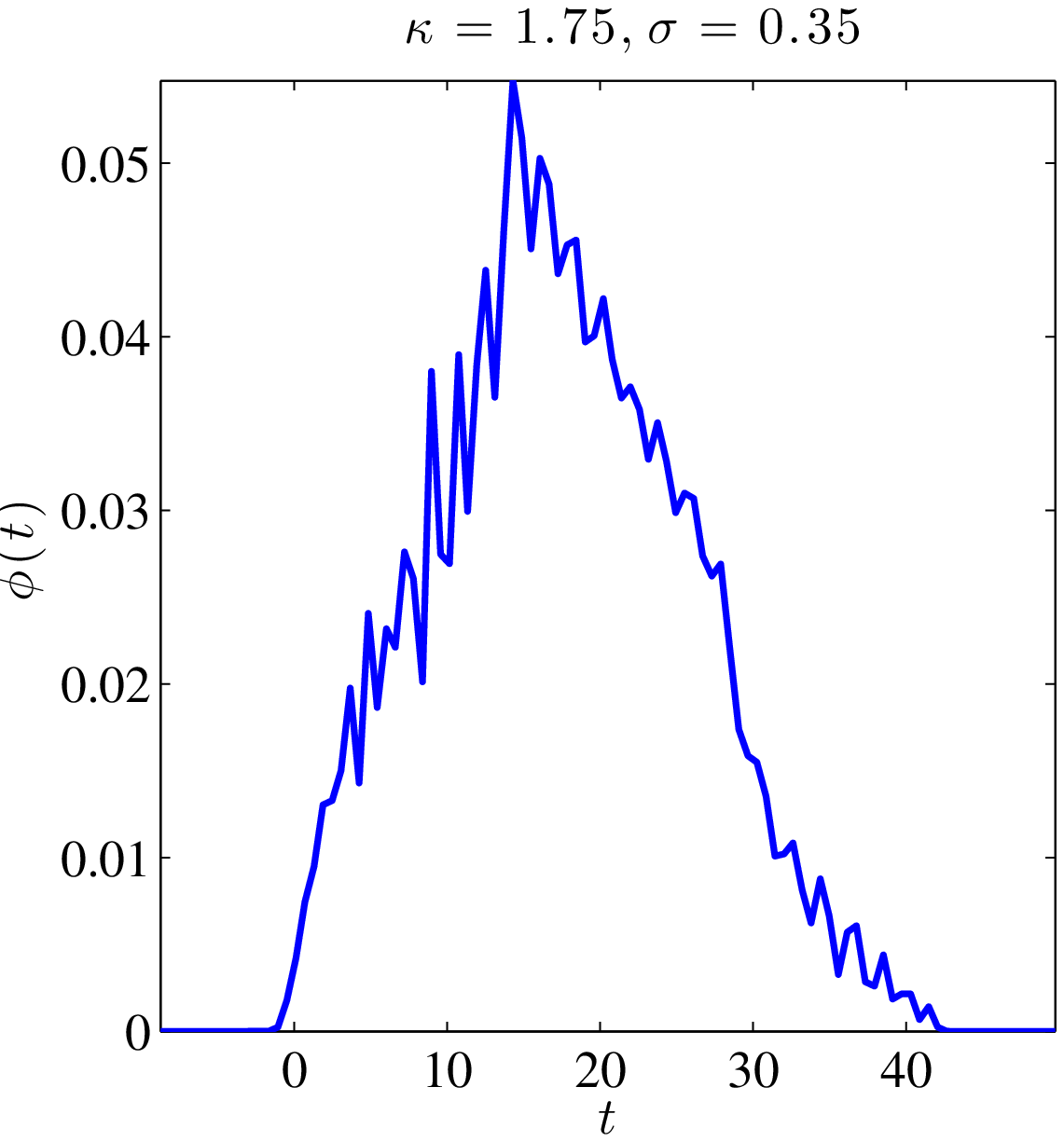}
\quad
\includegraphics[width=0.30\textwidth]{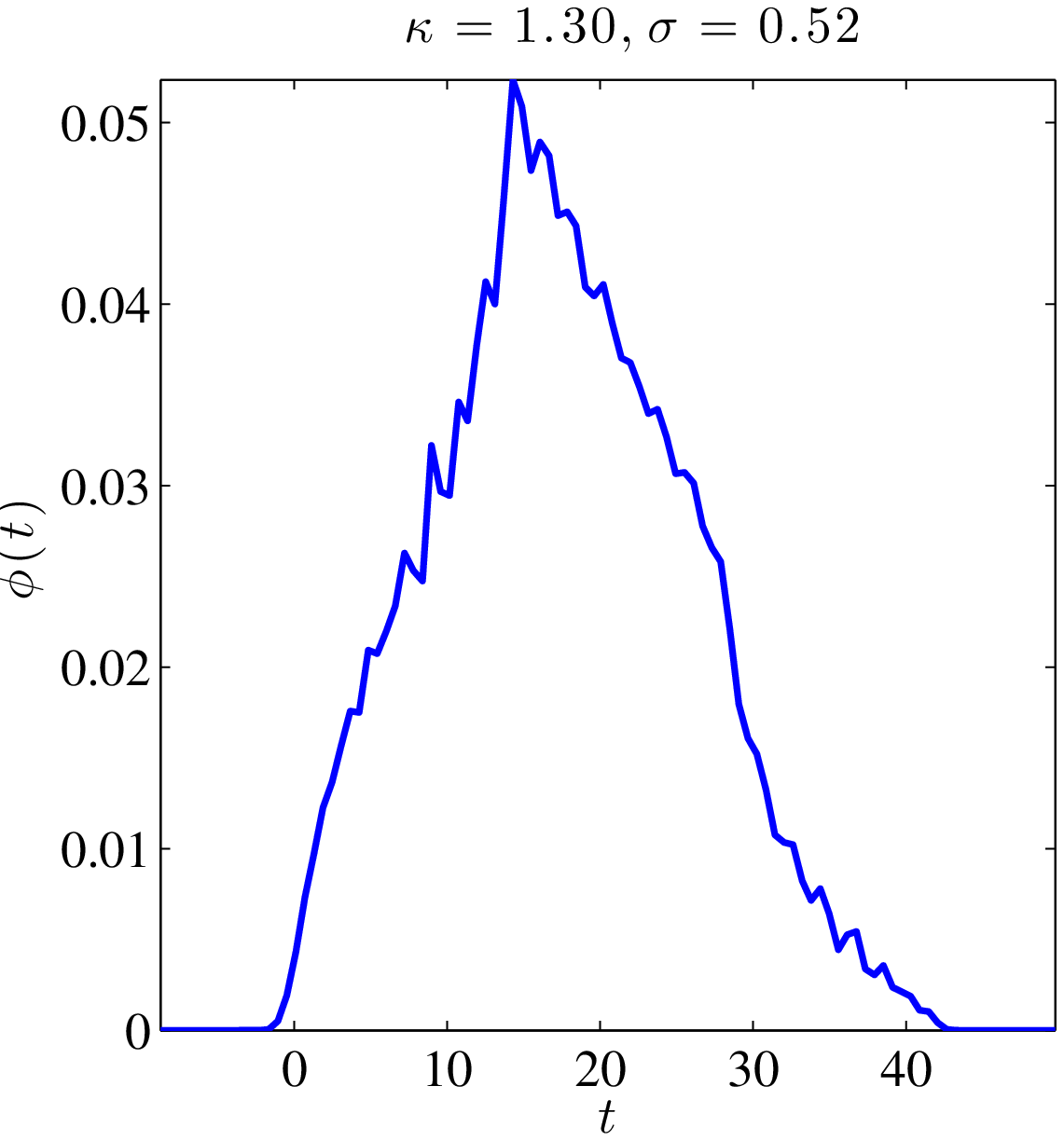}
 
\includegraphics[width=0.30\textwidth]{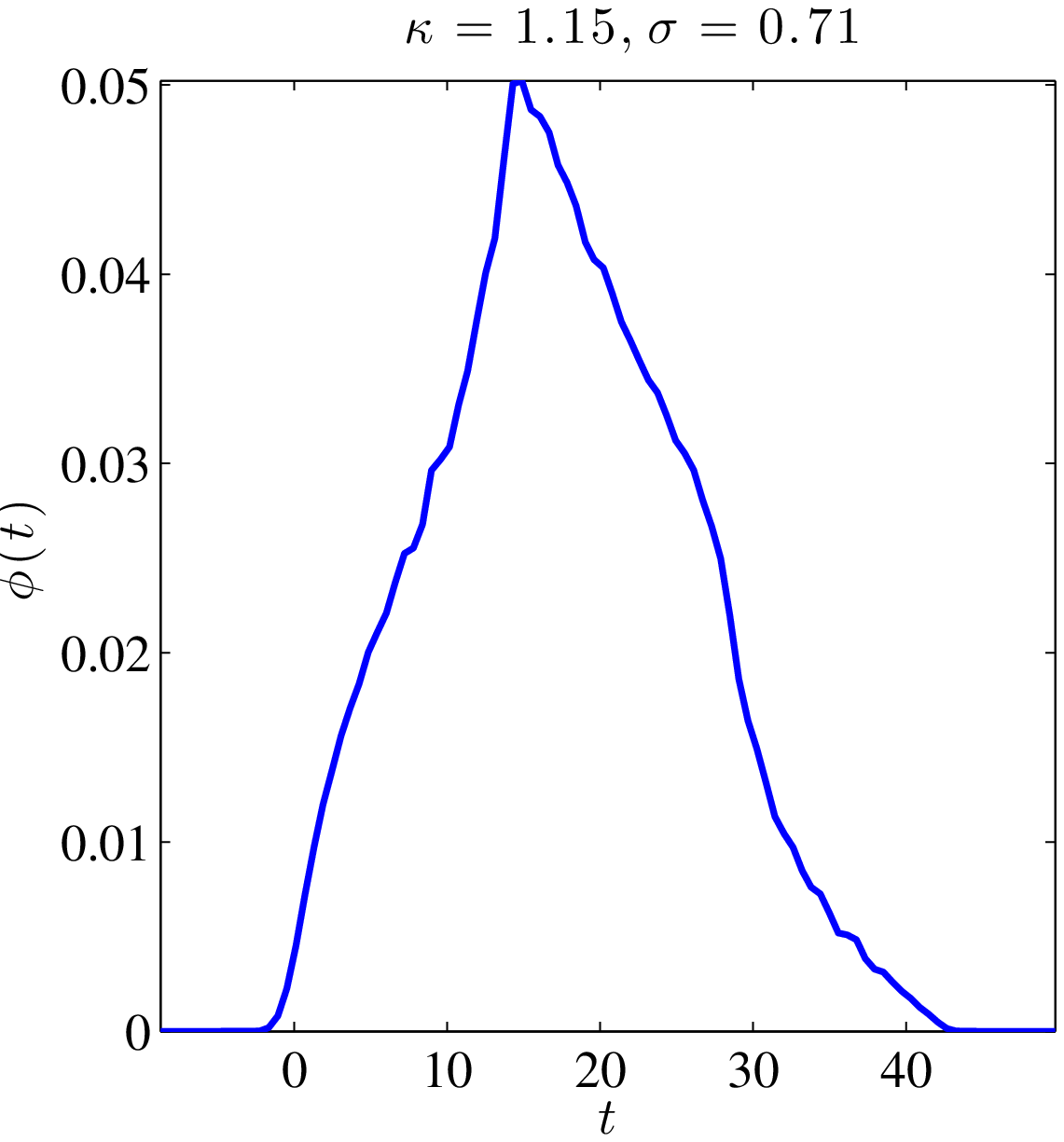}
\quad
\includegraphics[width=0.30\textwidth]{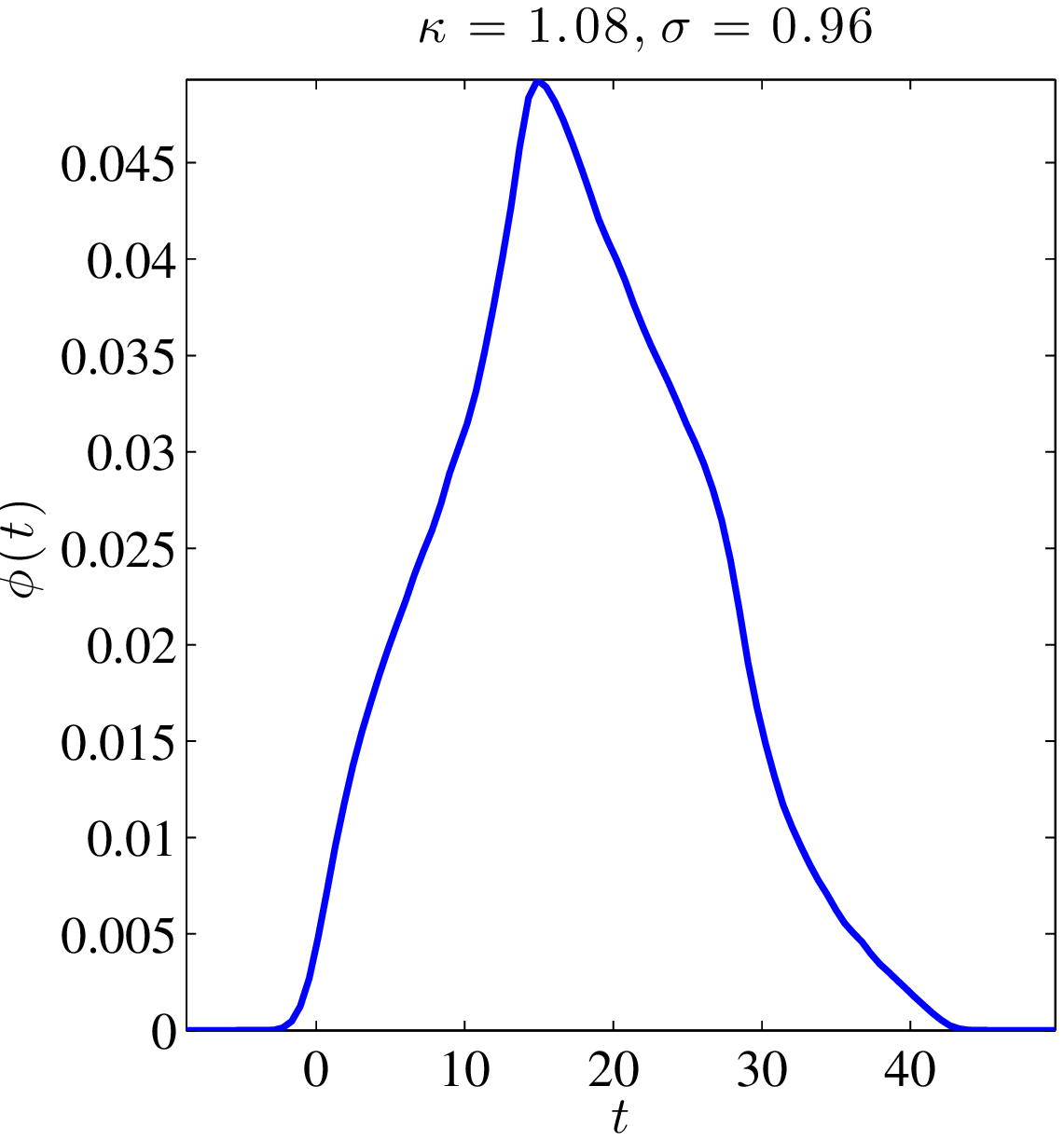}
	\end{center}
\caption{The eigenvalues, as well as various regularized
DOS $\phi_{\sigma}$ obtained by blurring the exact DOS (sum of
$\delta$-functions positioned at eigenvalues) of a matrix with Gaussians
of the form \eqref{eq:gauss}.}
\label{fig:si2Ex}
\end{figure}

The $\phi_{\sigma}(t)$ function defined in~\eqref{eq:blurphi} 
is one way to construct a regularized spectral density.
Again, the parameter $\sigma$ controls the resolution.
Larger values of $\sigma$ will lead to  smooth curves at the expense
of accuracy. Smaller values of $\sigma$ will lead to rough
curves that have peaks at the eigenvalues and zeros elsewhere.
This is illustrated in Figure~\ref{fig:si2Ex} where $\sigma$ takes
4 different values.  We can see that as $\sigma$ increases,
$\phi_{\sigma}$ becomes  smoother. When $\sigma = 0.96$,
which corresponds to a very smooth spectral density, we 
can still see the global profile of the eigenvalue distribution,
although local variation of the spectral density is mostly averaged out.

We remark that the optimal choice of $\sigma$, and therefore the
smoothness of the approximate DOS, is application dependent.  On 
the one hand, $\sigma$ should be chosen to be as large as possible so that the
regularized DOS $\phi_{\sigma}$ is easy to approximate numerically.
On the other hand, increasing $\sigma$ could cause an undesirable loss of detail
and  yield an erroneous result.  It is up to the user to select a value
of $\sigma$ that balances accuracy and efficiency: $\sigma$ should be
chosen to be small enough to reach the target accuracy, but not too
small which would require a large number of MATVECs to approximate the
DOS.

We should also note that ~\eqref{eq:blurphi} is not the only way
to regularize the DOS. Another choice is the Lorentzian function defined
as
\begin{equation}
\frac{\eta}{(t-\lambda)^2 + \eta^2} = -\mathrm{Im} \left(
\frac{1}{t-\lambda+i\eta}\right),
\label{eq:lorentz}
\end{equation}
where $\mathrm{Im}(z)$ denotes the imaginary part of a complex number $z$, and
$\eta$ is a small regularization constant that controls the width of the peak 
centered at $\lambda$.  As $\eta$ approaches zero,~\eqref{eq:lorentz} 
approaches a Dirac $\delta$-function centered at eigenvalues.
This approach is used in Haydock's method to be discussed in 
section~\ref{sec:haydock}.  We also examine the difference between the
regularization procedures in section~\ref{sec:examples} through
numerical experiments. 

\subsection{Non-negativity condition}

Since the DOS can be viewed as a probability distribution function, it
satisfies the non-negativity condition in the sense that
\begin{equation}
  \average{\phi,g} \ge 0,
  \label{eqn:nonnegative}
\end{equation}
for all non-negative function $g(t)\ge 0$ in the Schwartz space.  Not
all numerical methods described below satisfy the non-negativity
condition by construction.  We will see in section \ref{sec:examples}
that the failure of preserving the non-negativity condition can
possibly lead to large numerical error.

%
\section{Numerical methods for estimating spectral density} \label{sec:methods}
In this section, we review two classes of methods
for approximating the DOS of $A$.  We begin with the KPM, which can be 
viewed as a polynomial approximation to the DOS.  We show
that two other approaches that are derived from different view points
are equivalent to KPM. We then describe the second class of methods
which are based on the use of the familiar Lanczos partial tridiagonalization 
procedure~\cite{Lanczos1950}.
These methods use blurring (or regularization) techniques
to construct an approximate DOS from Ritz values.  They
differ in the type of blurring they utilize.  One of them, which we
will simply call the Lanczos method uses Gaussian blurring, whereas
the other method, which we call the Haydock's 
method~\cite{HaydockHeineKelly1972,epl2011,thermal2011,alloys2012}, 
uses a Lorentzian blurring.

A common characteristic of these methods is that they
all use a stochastic sampling and averaging technique to 
obtain an approximate DOS. The stochastic sampling and averaging
technique is based on the following key
result~\cite{Hutchinson1989,SilverRoder1994,AvronToledo2011}. 
\begin{theorem} \label{thm:stochastic}
Let $A$ be a real symmetric matrix of
dimension $n\times n$ with eigen-decomposition
$A = \sum_{j=1}^{n} \lambda_{j} u_{j} u_{j}^{T},$
and $u_{i}^{T}u_{j}=\delta_{ij},\quad i,j=1,\cdots,n$.  Here
$\delta_{ij}$ is the Kronecker $\delta$ symbol.  Let $v$ be a vector of
dimension $n$, and $v$ can be represented as the linear combination of
$\{u_{i}\}_{i=1}^{n}$ as
\begin{equation}
v = \sum_{j=1}^{n} \beta_{j} u_{j}. 
\label{eqn:linearcombination}
\end{equation}
If each component of $v$ is obtained independently from a normal
distribution with zero mean and unit standard deviation, i.e.
\begin{equation}
  \mathbb{E} [v]= 0, \quad \mathbb{E} [v v^{T}] = I,
  \label{eqn:stochasticcondition}
\end{equation}
then 
\begin{equation}
  \mathbb{E} [\beta_{i}\beta_{j}] = \delta_{ij}, \quad i,j=1,\cdots,n.
  \label{eqn:stochasticbeta}
\end{equation}
\end{theorem} 

The proof of Theorem~\ref{thm:stochastic} is straightforward.
The theorem suggest that the trace of a matrix function $f(A)$, 
which we need to compute in the KPM, for example, can be obtained
by simply averaging $v^{T} f(A)v$ for a number of randomly
generated vectors $v$ that satisfy the conditions given 
in~\eqref{eqn:stochasticcondition} because
\begin{equation}
\mathbb{E} [v^{T} f(A) v] = 
\mathbb{E} [\sum_{j=1}^{n} \beta_{j}^2 f(\lambda_{j})]= \sum_{j=1}^{n}
f(\lambda_{j}) = \mbox{Trace}[f(A)].
\label{eq:tracemean}
\end{equation}


\subsection{The Kernel Polynomial Method} \label{sec:kpm} 
The Kernel Polynomial Method (KPM) was proposed by Silver and 
R\"oder~\cite{SilverRoder1994} and Wang~\cite{Wang1994} in the mid-1990s to
calculate the DOS. See also
\cite{SilverRoder1997,SilverRoederVoterEtAl1996,DraboldSankey1993,ParkerZhuHuangEtAl1996}
among others where similar approaches were also used.

\subsubsection{Derivation of the origional KPM} \label{sec:kpmderive}
The KPM method constructs an approximation to the exact DOS of a matrix 
$A$ by formally expanding Dirac $\delta$-functions in terms of Chebyshev
polynomials $T_{k}(t)=\cos(k\arccos(t))$.
For simplicity, we assume that the eigenvalues are in the interval $[-1,
1]$.  As is the case for all methods which rely on Chebyshev expansions,
a change of variables must first be performed to map an interval that
contains $[ \lambda_{\min}, \lambda_{\max}]$ to $[-1, 1]$ if this assumption 
does not hold. Following
 the Silver-R\"oder paper~\cite{SilverRoder1994},  we include, for convenience,  
the inverse of the weight function into the spectral density function, 
\beeq{eq:tdos}
\hat \phi(t)  = 
\sqrt{1-t^2} \phi (t) =
\sqrt{1-t^2}  \times \frac{1}{n}
\sum_{j=1}^n  \delta(t - \lambda_j).
\en 
Then we expand the distribution $\hat{\phi}(t)$ as
\beeq{eq:exp1} 
\hat \phi(t)  = \sum_{k=0}^\infty  \mu_k T_k(t).
\en 
Eq.~\eqref{eq:exp1} should be understood in the sense of
distributions, i.e. for any test function $g\in \mathcal{S}$,
\[
\int_{-1}^{1}\hat{\phi}(t) g(t) \ dt = \int_{-1}^{1}
\sum_{k=0}^{\infty} \mu_{k} T_{k}(t) g(t) \ dt.
\]
The
same notation applies to the expansion of the DOS using other methods in the
following discussion.
By means of a formal moment matching procedure, the expansion
coefficients $\mu_k$ are also defined by
\begin{eqnarray} 
\mu_k & = &
 \frac{2-\delta_{k0}}{\pi} \int_{-1}^1 \frac{1}{\sqrt{1 - t^2}} T_k(t) 
\hat \phi (t) dt 
\nonumber \\
      & = &\frac{2-\delta_{k0}}{\pi}  \int_{-1}^1 \frac{1}{\sqrt{1 - t^2}} T_k(t) 
\sqrt{1-t^2} \phi (t) dt 
\nonumber \\
      & = &\frac{2-\delta_{k0}}{n\pi}  \sum_{j=1}^n  T_k(\lambda_j) . 
\label{eq:expmu} 
\end{eqnarray}
Here   $\delta_{ij}$ is the Kronecker $\delta$ symbol so that 
 $2 - \delta_{k0} $ is equal to 1 when $k=0$ and to 2 otherwise.


Thus, apart from the scaling factor $(2-\delta_{k0})/(n\pi)$, 
 $\mu_k$ is  the trace of $T_k(A)$. It follows from 
Theorem~\ref{thm:stochastic} that 
\begin{equation}
\zeta_{k} = \frac{1}{\nvec}  
\sum_{l=1}^{\nvec} \left(v^{(l)}_{0}\right)^T  T_k (A)
v^{(l)}_{0}
  \label{eqn:nuk}
\end{equation}
is a good estimation of the trace of $T_{k}(A)$,
for a set of randomly generated vectors $v^{(1)}_{0}$, $v^{(2)}_{0}$, ...,
 $v^{(\nvec)}_{0}$ that satisfy the conditions given by 
~\eqref{eqn:stochasticcondition}.
Here the subscript 0 is added to indicate that the vectors have not been
multiplied by the matrix $A$.
Then $\mu_k$ can be estimated by
\beeq{eq:muEst}
\mu_k \approx \frac{2-\delta_{k0}}{n\pi}  \zeta_{k}, 
\en

Now we consider the computation of each term $\left(v^{(l)}_{0}\right)^T
T_k (A) v^{(l)}_{0}$.  For simplicity we drop the superscript $l$ and
denote by $v_{0}\equiv v^{(l)}_{0}$. The 3-term recurrence of the
Chebyshev polynomial is exploited to compute $T_k(A)v_{0}$:
\[
T_{k+1} (A) v_{0} = 2 A T_k(A) v_{0} - T_{k-1}(A) v_{0}.
\]
So if we let $v_{k} \equiv T_k(A) v_{0}$, we have
\eq{eq:3term}
v_{k+1} = 2 A v_{k} - v_{k-1}.  
\en
Once the scalars $\{\mu_k\}$ are determined, we would in theory get the 
expansion for 
$ \phi(t) = \frac{1}{\sqrt{1 - t^2}} \hat{\phi} (t) $. In practice,
$\mu_{k}$ decays to $0$ as $k\to \infty$, and
the approximate density of states will be limited to Chebyshev polynomials of
degree $M$. So $\phi$ is approximated by:
\beeq{eq:tdosFinal} 
\wt{\phi}_M (t) = 
\frac{1}{\sqrt{1-t^2}} \sum_{k=0}^M \mu_k T_k (t) .
\eneq 

For a general matrix $A$ whose eigenvalues are not necessarily in the interval
$[-1,1]$, a linear transformation is first applied to $A$ to bring
its eigenvalues to the desired interval. Specifically, we will apply the 
method to the matrix
\[
B = \frac{A - c I}{d},
\]
where
\begin{equation}
c = \frac{\lambda_{lb} + \lambda_{ub}}{2} \ , 
\quad 
d = \frac{\lambda_{ub} - \lambda_{lb}}{2},
	\label{eqn:cd}
\end{equation}
and $\lambda_{lb}$, $\lambda_{ub}$ are lower and upper bounds of the 
smallest and largest eigenvalues $\lambda_{\min}$ and $\lambda_{\max}$ 
of $A$ respectively.

It is important to ensure that the eigenvalues of $B$ are within the
interval $[-1, 1]$. Otherwise the magnitude of the Chebyshev polynomial,
hence the product of $T_k(B)$ and $v_0$ computed through a three-term
recurrence, will grow exponentially with $k$. 

There are a number of ways~\cite{vandetal2000,ZhouSaadTiagoEtAl2006,LiZhou2011} to 
obtain good lower and upper bounds 
$\lambda_{lb}$ and $\lambda_{ub}$ of the spectrum of $A$.
For example, we can set 
$\lambda_{ub}$ to 
$\theta_k + \| (A - \theta_k I) u_k)\|$, 
and
$\lambda_{lb}$ to 
$\theta_1 - \| (A - \theta_1 I) u_1)\|$, 
where $\theta_1$ (resp. $\theta_k$) is the
algebraically smallest (resp. largest) Ritz value obtained from an $k$-step Lanczos iteration
and $u_1$ (resp. $u_k$) is the associated normalized Ritz vector.
Note that these residual norms are inexpensive to compute
since  $\| (A-\theta_j I)u_j\| $ can be easily expressed from 
the bottom entry of $z_j$,  the unit norm eigenvector of the $k \times k$ 
tridiagonal matrix obtained from the Lanczos process. For details,
see Parlett~\cite[sec. 13.2]{Parlett1998}.
We should point out that $\lambda_{lb}$ and $\lambda_{ub}$
do not have to be very accurate approximations to $\lambda_{\min}$ and
 $\lambda_{\max}$ respectively. It is demonstrated in~\cite{LiZhou2011} that 
tight bounds can be obtained from a 20-step Lanczos iteration for matrices
of dimension larger than 100,000.

\begin{algorithm}
\begin{small}
\begin{center}
  \begin{minipage}{5in}
\begin{tabular}{p{0.5in}p{4.5in}}
{\bf Input}:  &  \begin{minipage}[t]{4.0in}
                 Real symmetric matrix $A$ with eigenvalues between
                 $[-1,1]$. A set of points $\{t_{i}\}$
                 at which DOS is to be evaluated, the degree $M$ of
                 the expansion polynomial.
                 \end{minipage}\\
{\bf Output}: &  Approximate DOS $\{\wt{\phi}_{M}(t_i)\}$. \\
\end{tabular}
\begin{algorithmic}[1]
\STATE Set $\zeta_k = 0 $ for $k=0,\cdots, M$;
\FOR {$l=1:\nvec$}
\STATE Select a new random vector $v_0^{(l)}$;
\FOR {$k=0:M$}
\STATE Compute $\zeta_{k} \leftarrow \zeta_{k} + \left(v_0^{(l)}\right)^T
v_{k}^{(l)}$;
\STATE Compute $v_{k+1}^{(l)}$ via the three-term recurrence
$v_{k+1}^{(l)} =
2 A v_k^{(l)} - v_{k-1}^{(l)}$ \quad (for $k=0$, $v_1^{(l)} = A
v_0^{(l)}$);
\ENDFOR
\ENDFOR
\STATE Set $\zeta_k \leftarrow \zeta_k / \nvec, \mu_{k} \leftarrow \frac{2-\delta_{k0}}{n\pi}  \zeta_{k}$ for 
$k = 0,1,...,M$;
\STATE Evaluate $\{\wt{\phi}_{M}(t_i)\}$ using $\{\mu_{k}\}$ and
Eq.~\eqref{eq:tdosFinal};
\end{algorithmic}
\end{minipage}
\end{center}
\end{small}
\caption{The Kernel Polynomial Method.}
\label{alg:kpm}
\end{algorithm}

\begin{figure} 
	\begin{center}
    \includegraphics[width=0.35\textwidth,height=0.35\textwidth]{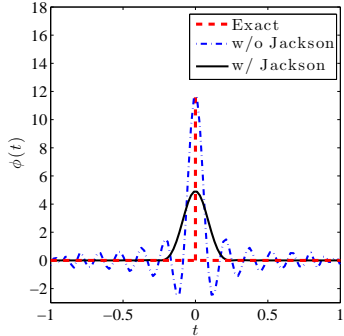}\quad
	\end{center}
\caption{Chebyshev expansion with and without the Jackson damping for Dirac
$\delta$-function $\delta(t)$. The Chebyshev polynomial degree is set to
$40$.}\label{fig:delCheb40}
\end{figure}

Because KPM can be viewed as a way to approximate a series of
$\delta$-functions which are highly discontinuous,
Gibbs oscillation~\cite{JeffreysJeffreys1999} can be 
observed near the peaks of the spectral density. 
Fig.~\ref{fig:delCheb40} shows an approximation to 
$\delta(t)$ by a Chebyshev polynomial expansion of the form
~\ref{eq:tdosFinal} with $M=40$.
It can been that the approximation oscillates around 0 away from $t=0$.
As a result, it does not preserve the nonnegativity of $\delta(t)$.  
The Gibbs oscillation can be reduced
by a technique called the Jackson damping, which modifies the coefficient of
Chebyshev expansion.  The details of Jackson damping is given in
Appendix \ref{sec:app_kpm}, and Fig.~\ref{fig:delCheb40} shows that the
Jackson damping indeed reduces the amount of Gibbs oscillation significantly. 
However, Jackson damping tends to oversmooth the approximate DOS and create 
an approximation to $\delta(t)$ that has a wider spread.  We will discuss
this again with numerical results in section~\ref{sec:examples}.


We should also point out that it is possible to replace Chebyshev
polynomials in KPM by other orthogonal polynomials such as the 
Legendre polynomials.  We will call the variant that uses Legendre 
polynomials to expand the spectral density by KPML.  

We also note that the cost for constructing KPM can be reduced by 
techniques presented in the  Appendix \ref{sec:app_kpm}.
Because KPM provides a finite polynomial expansion,  
the approximate DOS $\tphi$ can be evaluated at any arbitrary 
point $t$ once the $\mu_k's$ have been determined.

\subsubsection{The spectroscopic view of KPM} 

In his 1956  book titled ``Applied Analysis'' 
\cite{Lanczos1988}, Lanczos described a method for computing spectra of
real symmetric matrices, which he
termed ``spectroscopic''. This approach, which also
relies heavily on Chebyshev polynomials, is rather
unusual in that it assimilates the spectrum of a matrix
to a collection of frequencies and the goal is to detect these
frequencies by Fourier analysis. Because  it is not 
competitive with modern methods for computing eigenvalues,
this technique has lost its appeal.
However,  we will show below that the spectroscopic 
approach is well suited for computing approximate spectral densities, 
and it is closely connected to the KPM.

Let us assume that the eigenvalues of the matrix $A$ are in [-1,1].
The spectroscopy approach takes samples of a function of the form
\begin{equation}
f(t) = \sum_{j=1}^n \beta_j^2 \cos(\theta_j t),
\label{eq:lanprobe}
\end{equation}
where $\theta_{j}$'s are related to eigenvalues of $A$ according to
$\theta_j = \cos\lambda_j$, at $t=0, 1, 2,...,M$. Then one can
take the Fourier transform of $f(t)$ to reveal the spectral
density of $A$. If a sufficient number of samples are taken,
then the Fourier transform of the sampled function should have 
peaks near $\cos \lambda_j$, $j=1,2,...,n$, and an approximate spectral
density can be obtained.

Because $\lambda_j$'s are not known,~\eqref{eq:lanprobe} cannot be
evaluated directly. However, $M+1$ uniform samples of $f(t)$, i.e. $f(0)$, $f(1)$,
..., $f(M)$ can be obtained from the average of
\eq{eq:fseq}
v_0^Tv_0, \: v_0^T T_1(A) v_0, \: ..., \: v_0^T T_M(A) v_0,
\en
where $T_k(t)$ is the same $k$th degree Chebyshev polynomial of the first kind 
we used in the previous section, and $v_0$ is a random starting vector.

Taking a discrete cosine transform of~\eqref{eq:fseq} yields
\eq{eq:sp3} 
F(p) = \half \left(f(0) + (-1)^{p}f(M)\right) + \sum_{k=1}^{M-1} f(k) 
\cos\frac{k p \pi}{ M },\quad p=0, \cdots, M. 
\en 
Note that, as is customary, the end values are halved to account for the
discontinuity of the data at the interval boundaries. 
An approximation to the spectral density $\phi(t)$ can be obtained
from $F(p)$, $p = 0, 1, ..., M$ through an interpolation procedure.


We now show that the spectroscopic approach is closely
connected to KPM.  This connection can be seen by noticing that,
the coefficient $\zeta_{k}$ in Eq.~\eqref{eqn:nuk}
essentially gives an estimate of the following transform of the spectral
density $\phi(t)$ 
\begin{equation}
f(s) = \int \cos(s \arccos t) \phi(t) dt,
\label{eq:ftransform}
\end{equation}
evaluated at an integer $k$.

Since $t\in[-1,1]$, we can rewrite ~\eqref{eq:ftransform} as the continuous cosine transform
of a related function by introducing an auxiliary variable $\xi=\arccos t$, 
and define
\[
\psi(\xi) = \phi(\cos \xi) \sin \xi.
\]
It is then easy to verify that ~\eqref{eq:ftransform} can be written as
\[
f(s) = \int_{0}^{\infty} \cos (s\xi) \psi(\xi) \ d\xi.
\]
Thus $f(s)$ can indeed by obtained by performing a cosine transform
of $\psi(\xi)$.

If $f(s)$ is given, we can obtain $\psi(\xi)$ via the inverse cosine transform 
\begin{equation}
\psi(\xi) = \frac{2}{\pi} \int_{0}^{\infty} \cos(s\xi) f(s) \ ds.
\label{eq:finvpsi}
\end{equation}
Substituting $\xi=\arccos t$ back into \eqref{eq:finvpsi} yields
\begin{equation}
  \phi(t) = \frac{2}{\pi\sqrt{1-t^2}}\int_{0}^{\infty} \cos(s\arccos t) f(s) \ ds.
  \label{eqn:phifs}
\end{equation}

However, since we can only compute $f(s)$ for $s=k,k\in \mathbb{N}$ through 
estimating the trace of $T_{k}(A)$ by a stochastic averaging technique
discussed in section~\ref{sec:kpmderive}, the integration in \eqref{eqn:phifs} 
can only be performed numerically using, for example, a composite trapezoidal rule:
\begin{equation}
  \phi(t) \approx \frac{2}{n\pi\sqrt{1-t^2}} \left[ \frac12 f(0) + \sum_{k=1} f(k)
  T_{k}(t) + \frac12 f(M) T_{M}(t) \right]
  \label{eqn:phifs2}
\end{equation}
Comparing Eq.~\eqref{eqn:phifs2} with Eq.~\eqref{eq:exp1}, we find that
~\eqref{eqn:phifs2} is exactly the KPM expansion, except that the coefficient for
$T_{M}(t)$ is multiplied by a factor $\frac12$.  Therefore  the spectroscopic
method and KPM are essentially equivalent.

\subsubsection{The Delta-Gauss-Legendre expansion approach}\label{sec:dgl}
In some sense, the spectroscopy method discussed in the previous section 
samples the ``reciprocal" space of the spectrum by computing
$v_0^T T_j(A)v_0$, where $T_j(t)$ is the $j$th degree Chebyshev
polynomial of the first kind, for a number of different 
randomly generated vectors $v_0$ and
$j=0,1,...,M$. It uses the discrete cosine transform to reveal the 
spectral density.  In this section, we examine another way to sample 
or to probe the spectrum of $A$ at an arbitrary point
$t_i\in [-1,1]$ directly by computing 
$\{v_0^T p_{M_i}(A) v_0\}$, 
where $p_{M_i}(t)$ is an $M_i$th degree polynomial of the form
\eq{eq:DelCheb}
p_{M_i}(t) \equiv \sum_{k_i=0}^{M_i} \mu_{k_i}(t_i)
T_{k_i}(t).
\en 
The expansion coefficient $\mu_{k_i}(t_i)$ in the above expression 
is chosen, for each $t_i$, to be
\eq{eq:DelChebCoef}
\mu_{k_i}(t_i) 
= \frac{2-\delta_{k0}}{\pi} \int_{-1}^1 \frac{1}{\sqrt{1 - t^2}} T_{k_i}(t) 
 \delta (t - t_i) dt  
= 
 \frac{2-\delta_{k0}}{\pi} \frac{T_{k_i}(t_i)}{\sqrt{1 - t_i^2}}.
\en

The polynomial $p_{M_i}(t)$ defined in~\eqref{eq:DelCheb} can be viewed 
as a polynomial approximation to the $\delta$-function $\delta(t-t_i)$, 
which can be regarded as a spectral probe placed at $t_i$.  

The reason why $v_0^Tp_{M_i}(A) v_0$ can be regarded as a sample 
of the spectral density at $t_i$ can be explained as follows. 
The presence of an eigenvalue at $t_i$ can be
detected by integrating $\delta(t-t_i)$ over the entire spectrum 
of $A$ with respect to a spectral point measure defined at eigenvalues
only. 
The integral returns $+\infty$ if $t_i$ is an eigenvalue of 
$A$ and 0 otherwise. However, in practice, this integration cannot
be performed without knowing the eigenvalues of $A$ in advance.

A practical probing scheme can be devised by replacing $\delta(t-t_i)$ 
with a polynomial approximation such as the one given in~\eqref{eq:DelCheb}, 
and performing integration, which amounts to evaluating the
trace of $p_{M_i}(A)$.  This trace evaluation can be done by using the 
same stochastic approach we 
introduced earlier for the KPM. We call this technique
the {\em Delta-Chebyshev} method.

If $v_0$ is some random vector normalized to have $\|v_0\|=1$, 
whose expansion in the eigenbasis $\{u_j\}$ is given by 
\begin{equation}
v_0 = \sum_{j=1}^n \beta_j u_j,
\label{eq:v0}
\end{equation}
then $v_{M_{i}}\equiv p_{M_i}(A) v_0$ will have the expansion
\eq{eq:dsp} 
v_{M_i} = \sum_{j=1}^n \beta_j p_{M_i}(\lambda_j)  u_j. 
\en
Taking the inner product between $v_0$ and $v_{M_i}$ yields
\eq{eq:DelDens} 
 (v_{M_i})^T v_0 
= \sum_{j=1}^n \beta_j^2 p_{M_i}(\lambda_j).
\en 
Since $\sum_{j=1}^n \beta_j^2 = 1$, \eqref{eq:DelDens} can be 
viewed as an integral of $p_{M_i}(t)$ associated with
a point measure $\{\beta_j^2\}$ defined at eigenvalues of the 
matrix $A$.  In the case when $\beta_j^2 = 1/n$ for all $j$, we can 
then simply rewrite the integral as $\mbox{Trace}[p_{M_i}(A)]/n$.  
As we have already shown in previous sections, such a trace can be 
approximated by choosing multiple random vectors $v_0$ that
satisfy the conditions \eqref{eqn:stochasticcondition} and 
averaging $v_0^T p_{M_i}(A)v_0$ for all these  vectors. 
The averaged value yields
$\tphi(t_i)$, which is the approximation to the 
spectral density $\phi(t)$ at an arbitrary sample point $t_i$.

As we already indicated in section~\ref{sec:dosFun},  
because $\delta(t-t_i)$ is not a proper function, constructing
a good polynomial approximation directly may be difficult.
A more plausible approach is to ``regularize" $\delta(t-t_i)$ first
by replacing it with a smooth function that has a peak at $t=t_i$,
and constructing a polynomial approximation to this smooth
function.

We choose the regularized $\delta$-function 
to be the Gaussian $g_{\sigma}(t-t_i)$, where 
$g_{\sigma}$ is defined in \eqref{eq:gauss},
and the standard deviation $\sigma$ controls the smoothness
or the amount of regularization of the function.

It is possible to expand $g_{\sigma}(t-t_i)$ in terms of Chebyshev 
polynomials. However, we found that it is easier to 
derive an expansion in terms of Legendre polynomials.
It can been shown (see Appendix~\ref{sec:app_kpm}) that
\eq{eq:DGLexp} 
\gsz(t-t_i) = \frac{1}{(2\pi \sigma^2)^{1/2}} 
\sum_{k=0}^\infty \left(k+ \frac{1}{2} \right) \gamma_k(t_i)  L_k(t) \,
\en
where $L_k(t)$ is the Legendre polynomial of degree $k$,
and the expansion coefficient $\gamma_k(t_i)$ is defined by
\eq{eq:DGLcoeff} 
\gamma_k (t_i) = \ \int_{-1}^{1} L_k(s) e^{- \half ((s-t_i)/\sigma) ^2} ds .
\en
It can be also shown (see Appendix~\ref{sec:app_kpm}) that
$\gamma_k(t_i)$ can be determined by a recursive procedure that does not
require explicitly to compute the integral in~\eqref{eq:DGLcoeff}.

If we take an approximation to $\gsz(t-t_i)$ to be the first
$M_i+1$ terms in the expansion~\eqref{eq:DGLexp}, i.e.,
\begin{equation}
  \tilde{\phi}_{M_i}(t) = \frac{1}{(2\pi \sigma^2)^{1/2}}
  \sum_{k=0}^{M_i} \left(k+ \frac{1}{2} \right) \gamma_k(t_i)  L_k(t),
  \label{eqn:dglexpansion}
\end{equation}
then a practical scheme for sampling the spectral density of $A$
can be devised by computing $v_0^T h_{M_i}(A)v_0$ for randomly
generated and normalized $v_0$'s and averaging these quantities.
Because this scheme is based on regularizing the $\delta$ function
with a Gaussian and expanding the Gaussian in Legendre polynomials,
we call this scheme a Delta-Gauss-Legendre (DGL) method and summarize
this scheme in Algorithm~\ref{alg:DGL}.



\begin{algorithm}
\begin{small}
\begin{center}
  \begin{minipage}{5in}
\begin{tabular}{p{0.5in}p{4.5in}}
{\bf Input}:  &  \begin{minipage}[t]{4.0in}
                 Real symmetric matrix $A$ with eigenvalues between
                 $[-1,1]$. A set of points $\{t_{i}\}$
								 at which the DOS is to be evaluated, and $M_{\max}$ is
								 the maximum degree employed for all the points.
                 \end{minipage}\\
{\bf Output}: &  Approximate DOS $\{\wt{\phi}_{M}(t_i)\}$. \\
\end{tabular}
\begin{algorithmic}[1]
\FOR {each $t_i$}
\STATE Compute and store the expansion coefficients
$\{\gamma_{k}(t_{i})\}_{k=0}^{M_{i}}$ using Eq.~\eqref{eq:DGLcoeff};
\ENDFOR
\STATE Set $\zeta_k = 0 $ for $k=0,\cdots, M_{\max}$;
\FOR {$l=1:\nvec$}
\STATE Select a new random vector $v_0^{(l)}$;
\FOR {$k=0:M_{\max}$}
\STATE Compute $\zeta_{k} \leftarrow \zeta_{k} + \left(v_0^{(l)}\right)^T
v_{k}^{(l)}$;
\STATE Compute $v_{k+1}^{(l)}$ via the three-term recurrence
$v_{k+1}^{(l)} =
\frac{2k+1}{k+1} A v_k^{(l)} - \frac{k}{k+1} v_{k-1}^{(l)}$ \quad (for $k=0$, $v_1^{(l)} = A
v_0^{(l)}$);
\ENDFOR
\ENDFOR
\STATE Set $\zeta_k \leftarrow \zeta_k / \nvec$ for all $k = 0,1,...,M_{\max}$;
\STATE Evaluate $\wt{\phi}_{M_{i}}(t_i)$ using
Eq.~\eqref{eqn:dglexpansion} with $\{\zeta_k\}$ and 
the stored $\{\gamma_{k}(t_{i})\}$;

\end{algorithmic}
\end{minipage}
\end{center}
\end{small}
\caption{Multi-point Delta-Gauss-Legendre expansion.}
\label{alg:DGL}
\end{algorithm}

Note that both the Delta-Chebyshev and the DGL methods 
compute $v_0^T p_{M_i} v_0$ at sampled point $t_i$ within
the spectrum of $A$.  This would have been an unacceptably expensive 
procedure if it were not for the fact that the same vector
sequence $\{T_{k}(A) v_0\}$ and $\{L_k(A)v_0\}$ for $k=0,1,...$, 
can be used for all points $t_i$ at the same time. They only 
need to be generated once. 

Although the Delta-Chebyshev method and KPM are derived from 
somewhat different principles, there is a close connection
between the two, which may not be entirely obvious. 
The key to recognizing this connection is to notice that
the average value of $v_0^T p_{M_i}(A)v_0$ can be viewed as 
an approximation to $\mbox{Trace}(p_{M_i}(A))$, which can be
written as
\begin{eqnarray}
	\mbox{Trace}(p_{M_i}(A)) &=& \frac{1}{n}\sum_{k_i=0}^{M_i} \mu_{k_i}(t_i) \sum_{j=1}^n T_{k_i}(\lambda_j) \nonumber \\
	&=& \frac{1}{n} \sum_{k_i = 0}^{M_i} \frac{2-\delta_{k_i,0}}{\pi}
\frac{T_{k_i}(t_i)}{\sqrt{1-t_i^2}} \mbox{Trace}(T_{k_i}(A)) \nonumber \\
&=& \sum_{k_i = 0}^{M_i} \left[\frac{2-\delta_{k_i,0}}{n\pi}\mbox{Trace}(T_{k_i}(A)) \right] \frac{T_{k_i}(t_i)}{\sqrt{1-t_i^2}}. \label{eq:kpmform}
\end{eqnarray}
Note that the coefficients within the square bracket
in~\eqref{eq:kpmform} are exactly the same coefficients as in KPM, which
appear in the expansion~\eqref{eq:exp1} of the function $\hat{\phi}(t) =
\sqrt{1-t^2}\phi(t)$. 
Therefore, when $M_i = M$ for all $i$, the Delta-Chebyshev expansion
method is identical to the KPM. Hence, the cost of this approach 
is the same as that of KPM if polynomials of the same degree and the same 
number of sampling vectors are used at each $t_i$.

When $M_i$ is allowed to vary with respect to $i$, there is a 
slight advantage of using the Delta-Chebyshev method in terms of flexibility. 
We can use polynomials of different degrees in different parts of the 
spectrum to obtain a more accurate approximation.  Note that in this 
situation, if $M_{\max}$ is the maximum degree employed for all the points, 
the number of MATVECs employed remains the same 
and equal to $M_{\max}$, since we will need to compute for each 
random vector $v_0$, the vectors $T_k (A) v_0$ for $k=0,\cdots,M_{\max}$ 
as these are needed by the points requiring the highest degree. 
However, some of the other calculations (inner products) required to 
obtain the spectral density can be avoided, though in most cases
applying $T_{k}(A)$ to $v_{0}$ dominate the computational cost in the
DOS calculation.   The computational cost of DGL is similar to that
of the Delta-Chebyshev method. Similarly, one can also show that DGL is
closely related to KPML, i.e., it is expansion of a regularized 
spectral density in terms of Legendre polynomials. However, we will 
omit the alternative derivation here.

The close connection between the Delta-Chebyshev method and KPM also suggests that
Gibbs oscillation can be observed in the approximate DOS 
produced by Delta-Chebyshev and DGL, especially when $\sigma$ is small. There is
no guarantee that the non-negativity of $\phi(t)$ can be preserved 
by DGL.

\subsection{The Lanczos Algorithm} \label{sec:lanczos}

Because finding a highly accurate DOS essentially amounts to computing
all eigenvalues of $A$, any method that can provide approximations to
the spectrum of $A$ can be used to construct an approximate DOS as well.
Since the Lanczos algorithm yields good approximation to extreme
eigenvalues, it is a good candidate for computing localized spectral
densities at least at both ends of the spectrum.  In this section, we
show that it is also possible to combine the Lanczos algorithm with
multiple randomly generated starting vectors to construct a good
approximation to the complete DOS.  

It should be noted that the spectral density $\phi$ as a probability
distribution is non-negative, i.e. $\average{\phi,g}\ge 0$ if $g\ge 0$
everywhere. This is an important property, but the KPM and its variants
as introduced in previous sections do not preserve the non-negativity of
the spectral density.  In contrast, the methods introduced in this
section, including the Lanczos method and the Haydock method do preserve
the non-negativity by construction.  This will become a clear advantage
for certain spectral densities as will be illustrated in
section~\ref{sec:examples} through numerical experiments.

\subsubsection{Constructing spectral density approximating from Ritz values and Gaussian blurring} 
For a given starting vector $v_0$, an $M$-step Lanczos procedure for a
real symmetric matrix $A$ can be succinctly described by 
\begin{equation}
  AV_M = V_M T_M + f_M e_M^T, \ \ V_M^TV_M = I, \ \ V_M^Tf_M = 0.
  \label{eq:lanfact}
\end{equation}
Here $T_M$ is an $M\times M$ tridiagonal matrix, $V_{M}$ is an
$n\times M$ matrix, and $I_{M}$ is an $M\times M$
identity matrix. It is well known that~\cite{GVL-book} the $k$-th column 
of $V_M$ can be expressed as
\[
V_M e_k = p_{k-1}(A) v_0, \quad k=1,\cdots,M. 
\]
where $\{p_k(t)\}$, $k = 0,1,2,..,M-1$ is a set of polynomials orthogonal
with respect to the weighted spectral distribution
$\phi_{v_0}(t)$ that takes the form
\begin{equation}
	\phi_{v_0}(t) = \sum_{j=1}^n \beta_{j}^2 \delta(t-\lambda_{j}),
	\label{eq:v0dist}
\end{equation}
where $\beta_j$'s are the expansion coefficients obtained from 
expanding $v_{0}$ in the eigenvector basis of $A$ as 
in Eq.~\eqref{eq:v0}.

It is also well known that these orthogonal polynomials can be generated by a 
three-term recurrence whose coefficients are defined by the matrix elements 
of $T_{M}$~\cite{gautschi67}.  If $(\theta_{k},y_{k})$, $k=0,1,2,\ldots,M$ are eigenpairs of the 
tridiagonal matrix $T_M$, and $\tau_{k}$ is the first
entry of $y_{k}$, then the following distribution function defined by
\begin{equation}
	\sum_{k=0}^M \tau_{k}^2 \delta(t-\theta_k),
\label{eq:landist}
\end{equation}
serves as an approximation to the weighted spectral density function
$\phi_{v_0}(t)$, in the sense that
\begin{equation}
\sum_{j=1}^n \beta_{j}^2 p_q(\lambda_j) = \sum_{k=0}^M \tau_{k}^2 
p_q(\theta_k), 
\label{eq:tauexpand}
\end{equation}
for all polynomials of degree $0 \leq q \leq 2M+1$.  The moment matching
property described by~\eqref{eq:tauexpand} is well known \cite{gautschi81}, 
and is related to the Gaussian quadrature rules
\cite{gautschi68,golubwelsch,GVL-book}.

Since in most cases, we are interested in the standard spectral 
density defined by~\eqref{eq:DOS0}, we would like to choose a starting 
vector $v_0$ such that $\beta_{j}^2$ is uniform. However, this is generally
not possible without knowing the eigenvectors $\{u_j\}$ of $A$ in advance.
To address this issue, we resort to the same stochastic approach we used
in previous sections. 

We repeat the Lanczos process with multiple randomly generated starting 
vectors $v^{(l)}_0$, $l = 1,2,\ldots,\nvec$, that satisfy the
conditions given by~\eqref{eqn:stochasticcondition}. 
It follows from \eqref{eq:tracemean} that 
\begin{equation}
		\frac{1}{\nvec n}\sum_{l=1}^{\nvec} \left( v^{(l)}_{0}\right)^{T} \delta(tI-A) v^{(l)}_{0}	
= 
\frac{1}{n}\sum_{j=1}^{n}
		\left(\frac{1}{\nvec}\sum_{l=1}^{\nvec} \left(\beta_{j}^{(l)}\right)^2\right)
		\delta(t-\lambda_{j})
	\label{eq:avgphi}
\end{equation}
is a good approximation to the 
standard spectral density $\phi(t)$ in Eq.~\eqref{eq:DOS0}. 
Since each distribution~\eqref{eq:landist} generated by the Lanczos procedure 
is a good approximation to~\eqref{eq:v0dist}, the average of~\eqref{eq:landist} over $l$, i.e.,
\begin{equation}
\wt{\phi}(t) = \frac{1}{\nvec} \sum_{l=1}^{\nvec}\left(
\frac{1}{n} \sum_{k=0}^{M}  \left(\tau_{k}^{(l)}\right)^2
\delta(t-\theta_{k}^{(l)}) \right)
\label{eq:lanapprox}
\end{equation}
should yield a good approximation to the standard spectral density~\eqref{eq:DOS0}.


Since \eqref{eq:lanapprox} has far fewer peaks than $\phi(t)$,
when $M$ is small, a direct comparison of \eqref{eq:lanapprox} with 
$\phi(t)$ is not very meaningful.  However, when $\phi(t)$ is 
regularized by replacing $\delta(t-\lambda_i)$ with $g_\sigma(t-\lambda_i)$,
we can replace $\delta(t-\theta_k^{(l)})$ in~\eqref{eq:lanapprox} with a
Gaussian centered at the $\theta_k^{(l)}$ to yield a regularized DOS 
approximation, i.e., we define the approximate DOS as
\begin{equation}
	\tphi_{\sigma}(t) = \frac{1}{\nvec}\sum_{l=1}^{\nvec}\left(\frac{1}{n}
	\sum_{k=0}^{M} \left(\tau_{k}^{(l)}\right)^2
	\gsz(t-\theta_{k}^{(l)})\right).
	\label{eqn:landos}
\end{equation}
This regularization is well justified because in the limit of $M = n$,
all Ritz values are the eigenvalues, and $\tphi_{\sigma}(t)$ is exactly the
same as the regularized DOS $\phi_{\sigma}(t)$ for the same $\sigma$.
We will refer to the method that constructs the DOS approximation from 
Ritz values obtained from an $M$-step Lanczos iteration as the Lanczos 
method in the following discussion.  

Because $g_{\sigma}(t) \ge 0$, the approximate DOS produced by
the Lanczos method is nonnegative. This is a desirable property 
not shared by KPM, the DGL or the spectroscopic method.

An alternative way to refine the Lanczos based DOS approximation from a
$M$-step Lanczos run is to first construct an approximate cumulative
spectral density or cumulative density of states (CDOS), which is 
a monotonically increasing function and then take the derivative of 
the CDOS through a finite difference procedure or other means.  
This technique is discussed in Appendix \ref{sec:app_lanczos}.

\subsubsection{Haydock's method} \label{sec:haydock}
As we indicated earlier, the use of Gaussians is not the only way
to regularize the spectral density. Another possibility is
to replace $\delta(t-\lambda_i)$ in \eqref{eq:DOS1} with a Lorentzian 
of the form \eqref{eq:lorentz} and centered at $\lambda_i$. 
The regularized DOS can be written as
\[
\phi_{\eta}(t) =\frac{1}{n\pi}\sum_{j=1}^n \frac{\eta}{(t-\lambda_j)^2+\eta^2}.
\]

Consequently, an alternative approximation to the spectral density 
can be obtained by simply replacing $\delta(t-\theta_k^{(k)})$ in
\eqref{eq:lanapprox} with a Lorentzian centered at $\theta_k^{(k)}$,
i.e.,
\[
\tphi_{\eta}(t) = \frac{1}{\nvec}\sum_{l=1}^{\nvec}\left[\frac{1}{n}
	\sum_{k=0}^{M} \left(\tau_{k}^{(l)}\right)^2
	\frac{\eta}{(t-\theta_k^{(l)})^2 + \eta^2}\right],
\]
where $\theta_{k}^{(l)}$ and $\tau_{k}^{(l)}$ are the same
Ritz values and weighting factors that appear in~\eqref{eqn:landos} 
and $\eta$ is an appropriately chosen constant that corresponds
to the resolution of the spectral density to be approximate.
This approximation was first suggested by Haydock, Heine and
Kelly~\cite{HaydockHeineKelly1972}.  We will refer to this approach
as Haydock's method.  

Haydock's original method does not require computing Ritz values
even though computing the eigenvalues of a small tridiagonal 
matrix is by no means costly nowadays.  The method makes use 
of the fact that
\[
\phi_{\eta}(t) =-\frac{1}{n\pi} \Im \sum_{j=1}^{n} \frac{1}{t-\lambda_{j}+i\eta} 
= \frac{1}{n\pi} \Im  \mbox{Trace}\left[(tI-A+i\eta I)^{-1} \right].
\]
Hence, once again, the task of approximating $\phi_{\eta}(t)$ reduces to 
that of approximating the trace of $(tI-A+i\eta I)^{-1}$, 
which can be obtained by
\begin{equation}
\frac{1}{\nvec}\sum_{l=1}^{\nvec} \left(v_0^{(l)}\right)^T(t I-A+i\eta I)^{-1}
	v_0^{(l)},
\label{eqn:trinva}
\end{equation}
for ${\nvec}$ randomly generated vectors $v_0^{(l)}$ that satisfy the 
conditions \eqref{eqn:stochasticcondition}.

Note that a direct calculation of~\eqref{eqn:trinva} requires solving
linear systems of the form $[A - (t + i\eta)I] z = v_0$ repeatedly for 
any point $t$ at which the spectral density is to be evaluated. 
This approach can be prohibitively expensive.
Haydock's approach approximates
$v_0^T(t I -A+i\eta I)^{-1} v_0$ for multiple $t$'s 
at the cost of performing a single Lanczos factorization and some
additional calculations that are much lower in complexity.

If $v_0$ is used as the starting vector of the Lanczos procedure,
then it follows from the shift-invariant property of the Lanczos algorithm
that
\begin{equation}
[A - (t + i\eta)I] V_M = V_M [T_M - (t + i\eta)I] + f e_{M+1}^T, 
\label{laninvfac}
\end{equation}
where $V_M$ and $T_M$ are the same orthonormal and tridiagonal matrices 
respectively that appear in~\eqref{eq:lanfact}.  After multiplying
~\eqref{laninvfac} from the left by $[A - (t + i\eta)I]^{-1}$, from 
the right by $[T_M - (t + i\eta)I]^{-1}$ and rearranging terms, we
obtain
\[
[A - (t + i\eta)I]^{-1} V_M = V_M [T_M - (t + i\eta)I]^{-1} 
-[A - (t + i\eta)I]^{-1}fe_{M+1}^T  [T_M - (t + i\eta)I]^{-1}.
\]
It follows that  
\begin{eqnarray*}
 v_0^T\left[A-(t+i\eta) I \right])^{-1} v_0 &=& 
e_1^T V_M^T [A - (t + i\eta)I]^{-1} V_M e_1 \\
&=& e_1^T [T_M - (t + i\eta)I]^{-1} e_1 + \xi,
\end{eqnarray*}
where $\xi = - \left(v_0^T [A - (t + i\eta)I]^{-1}f\right) \left(e_{M+1}^T  
[T_M - (t + i\eta)I]^{-1} e_1\right)$.  If $\xi$ is sufficiently small, computing 
$v_0^T(t I-A+i\eta I)^{-1} v_0$ 
reduces to computing the $(1,1)$-th entry of the inverse
of $ T_M - (t + i\eta)I$.  It is not difficult to show that this
entry is exactly the same as the expression given in \eqref{eqn:trinva} 
up to a constant scaling factor.

Because $T_M$ is tridiagonal 
with $\alpha_1, \alpha_2, ..., \alpha_{M}$ on the diagonal and 
$\beta_2, \beta_3, ..., \beta_{M}$ on the sub-diagonals and super-diagonals, 
$e_1^T(zI - T_{M})^{-1} e_1$ can be computed in a recursive 
fashion using a continued fraction formula
\begin{equation}
	e_1^T(zI - T_{M})^{-1} e_1 = \frac{1}{z - \alpha_{1} +
	\frac{\beta_{2}^2}{z - \alpha_{2} + \cdots}}.
\label{eqn:cfHaydock}
\end{equation}
This formula can be verified from the identity
\[	\left( z - T_{M} \right)^{-1}_{1,1} \equiv
\frac{\det (z I - T_{M} ) }{ \det (z - \hat{T}_{M})}
\]
where $\hat{T}_{M}$ is the trailing submatrix starting from the
$(2,2)$ entry of $T_{M}$ \nref{eqn:cfHaydock} and tridiagonal structure
of both $T_{M}$ and $\hat{T}_{M}$ matrices. 

It is also related to the generation of Sturm sequences which is used in
bisection methods for computing eigenvalues of tridiagonal matrices
\cite{Parlett1998}.
Although this is an elegant way to compute $e_1^T (T_M - z I)^{-1}e_1$,
its cost is not much lower than that of solving the linear system 
$(T_M - z I)w=e_1$  and taking its first entry. For most
problems, the cost for this procedure is small compared to that required 
to perform the Lanczos procedure to obtain $T_M$.

We should point out that the non-negativity of $\phi(t)$ is preserved by 
the Haydock method.  However, the Lorentzian function defined 
by~\eqref{eq:lorentz} decreases to 0 at a much slower rate compared to 
a Gaussian as $t$ moves away from $\lambda$. Hence, when a high resolution 
approximation is needed, we may need to choose a very small $\eta$ in order 
to produce an accurate approximation.

\section{Numerical results}\label{sec:examples}

In this section, we compare all methods discussed above for
approximating the spectral density of $A$ through numerical examples. 
For a given test problem, we first define the target resolution
of the DOS to be approximated by setting the parameter 
$\sigma$ in either~\eqref{eq:err2} or~\eqref{eqn:errormetric} 
or the parameter $\eta$ in \eqref{eq:lorentz}.
We use the metric defined in Eq.~\eqref{eqn:errormetric} to measure the 
approximation errors associated with KPM and its variants with the
exception of DGL. For
DGL, Lanczos and Haydock methods, we simply use 
the error metric~\eqref{eq:err2} with $p=\infty$. Since the 
spectroscopic method is equivalent to KPM, we do not
show any numerical results for the spectroscopic method.


\subsection{Modified Laplacian matrix}\label{sec:lapmatrix}

The first example is a modified 2D Laplacian 
operator with zero Dirichlet boundary condition defined on
the domain $[0, 30] \times [0,30]$. The operator 
is discretized using a five-point finite difference stencil with
$\Delta h = 1$.  The modification involves adding
a diagonal matrix, which can be regarded as a discretized potential 
function. The diagonal matrix is generated by adding two Gaussians, 
one centered at the point (4,5) of the domain and the other at the 
point (25,15).  The dimension of the matrix 
is $750$, which is relatively small.  We set the parameter 
$\sigma$ and $\eta$ to $0.35$.   For all calculations shown in this section, 
we use $\nvec=100$ random vectors whenever stochastic averaging is 
needed.  Each calculation is repeated $10$ times.  Each plotted
value is the mean value of the computed quantities produced from the $10$ runs, 
with the error bar indicating the standard deviation of the $10$ runs.

\begin{figure}[h]
  \begin{center}
    \includegraphics[width=0.4\textwidth]{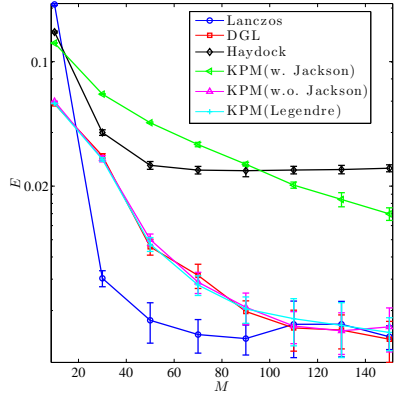}
  \end{center} 
  \caption{A comparison of approximation errors of all methods
           applied to the modified Laplacian matrix for different $M$
           values.}
  \label{fig:lap2d_all_blur}
\end{figure}

\begin{figure}[h]
  \begin{center}
    \includegraphics[width=0.4\textwidth]{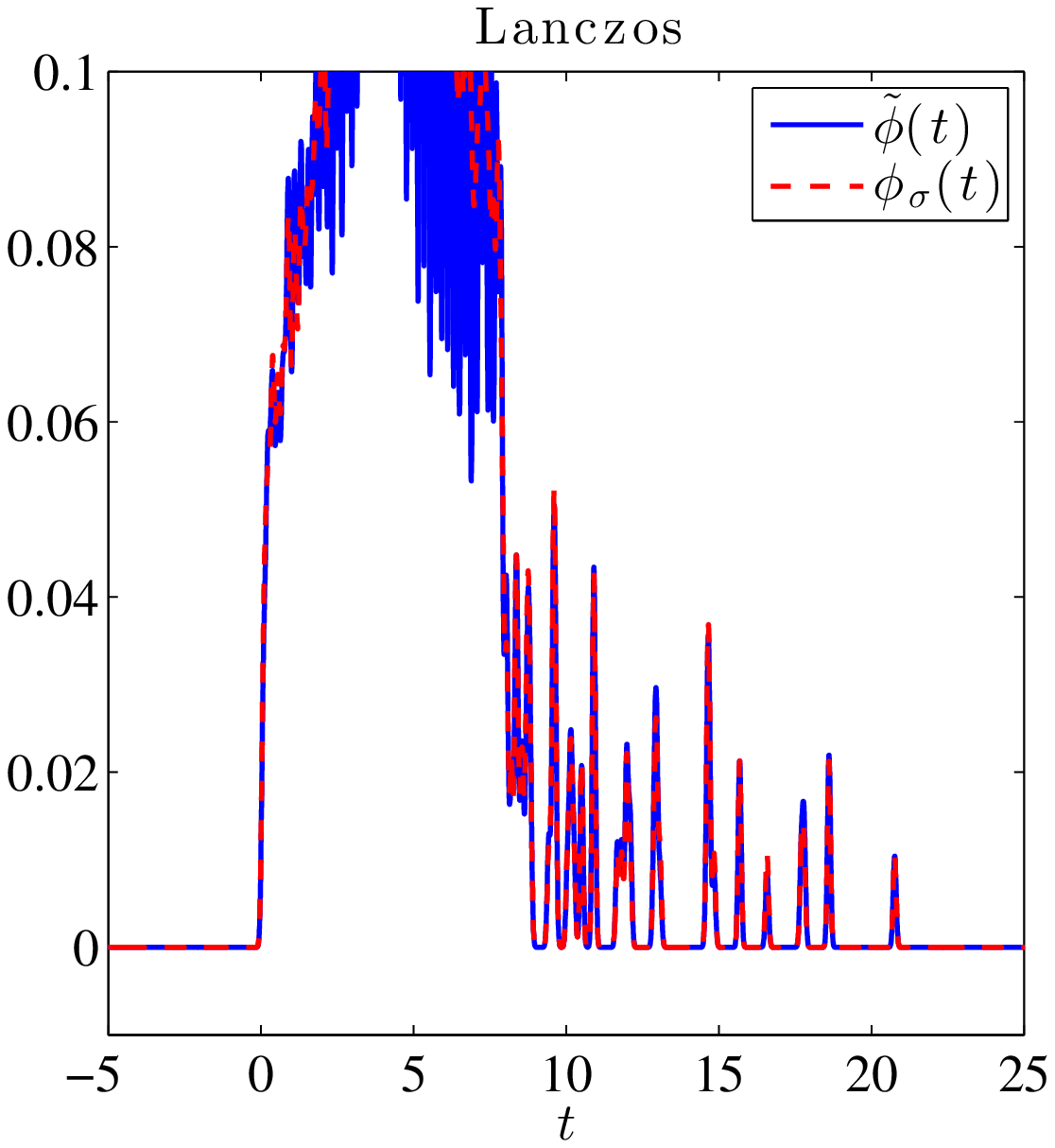}
    \includegraphics[width=0.4\textwidth]{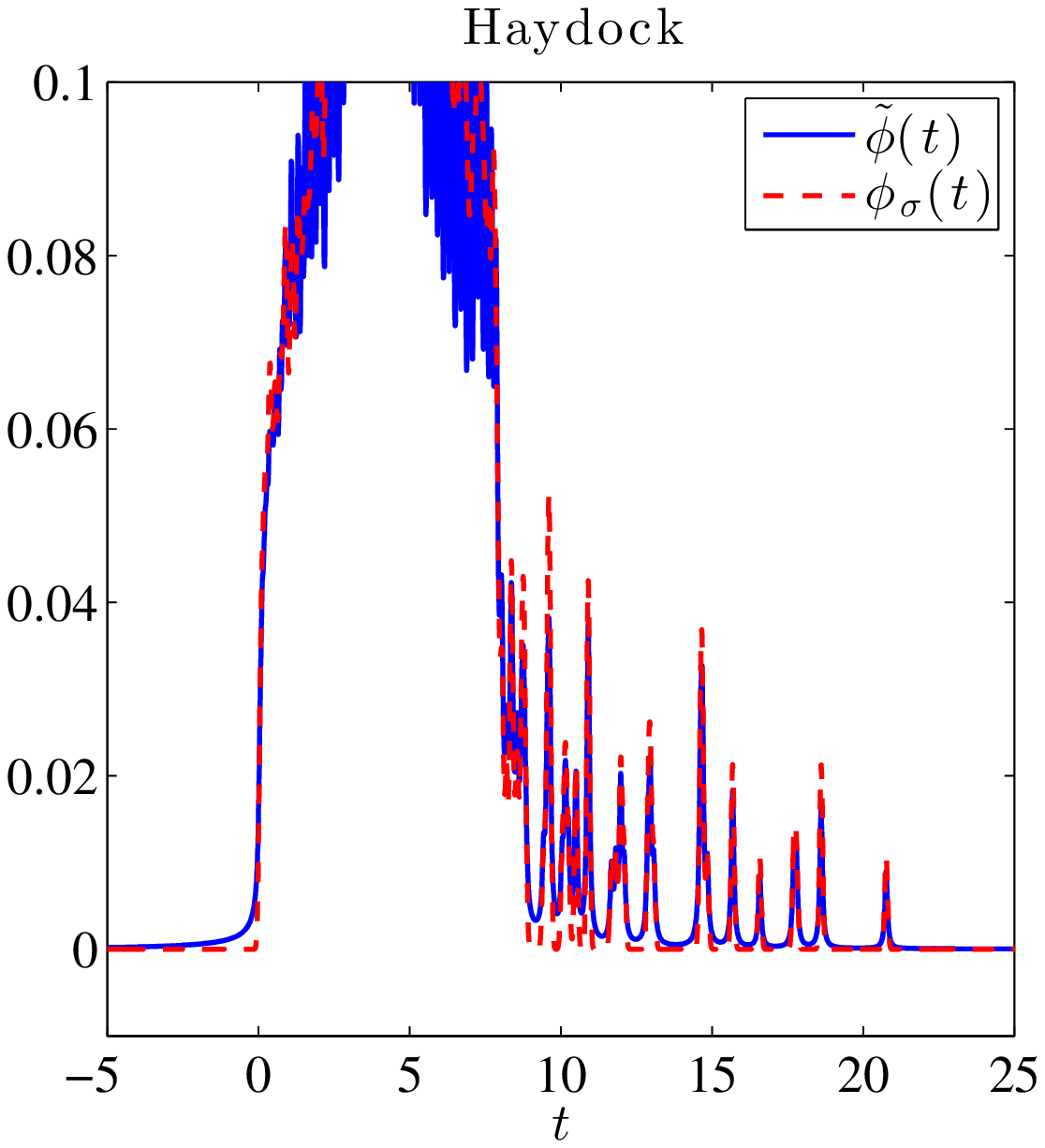} 
    
    \includegraphics[width=0.4\textwidth]{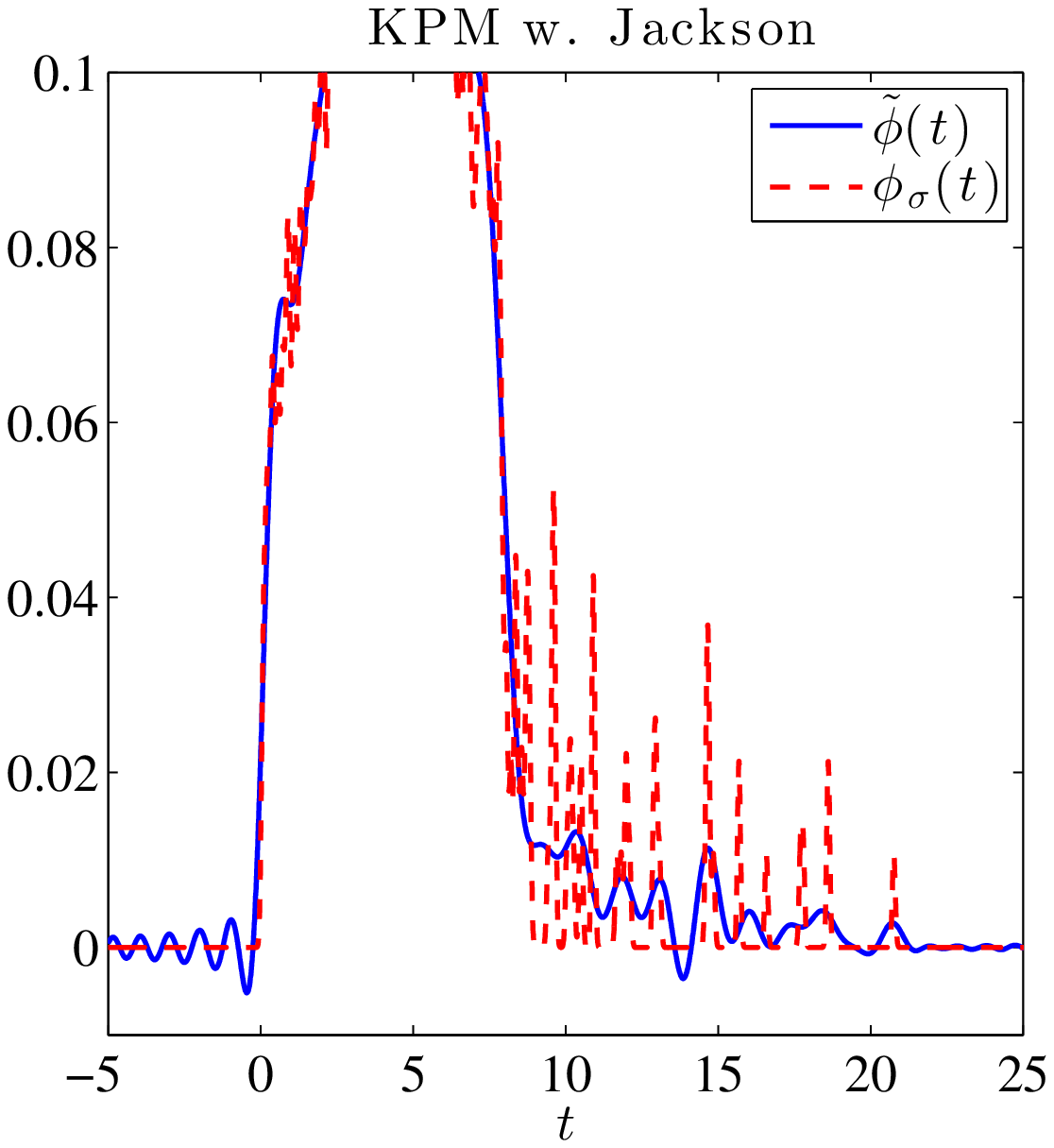}
    \includegraphics[width=0.4\textwidth]{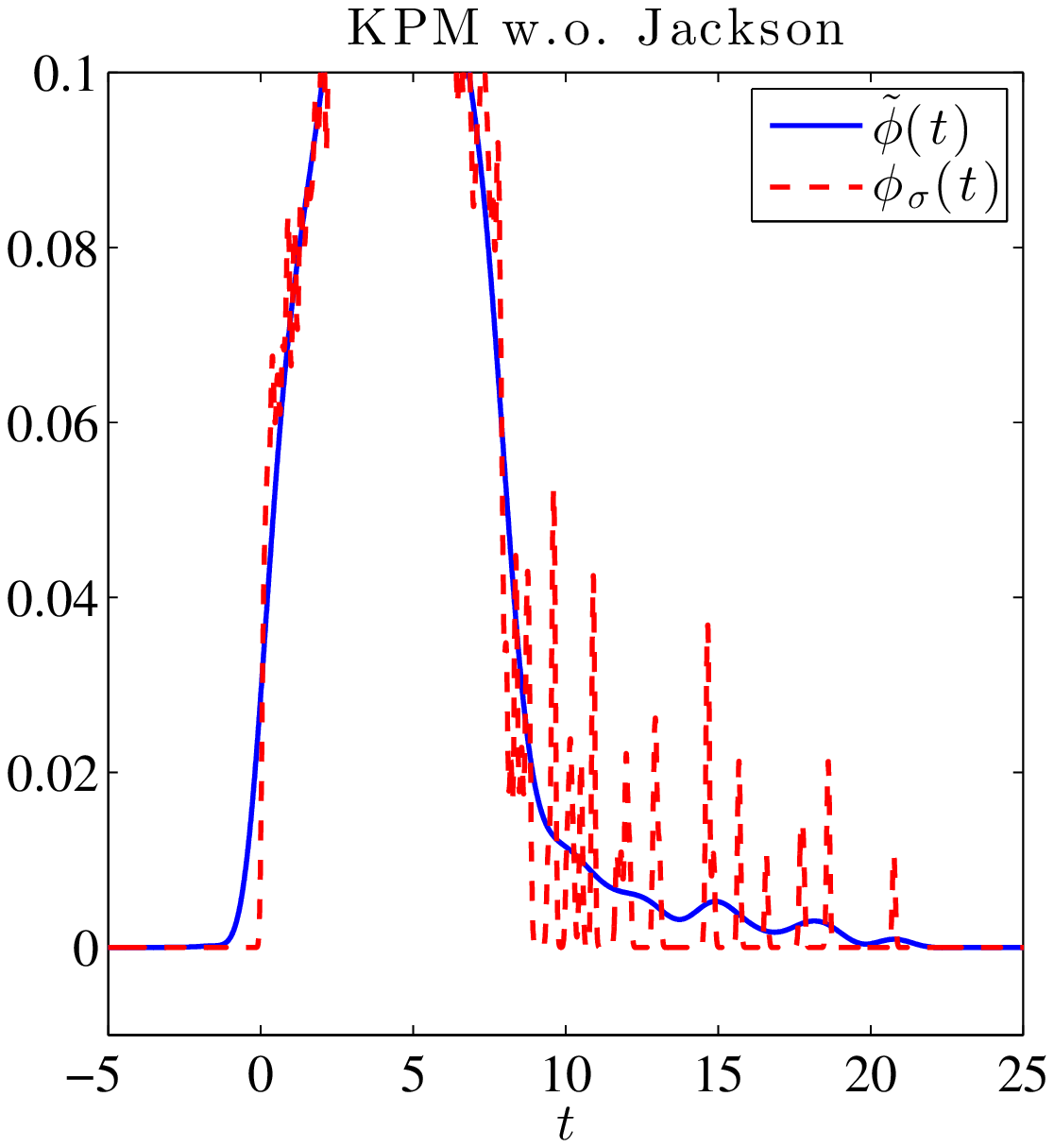}
  \end{center}
  \caption{Comparing the regularized DOS (with $\sigma=0.05$) with
           the approximate DOS produced by (a) the Lanczos method
           (b) the Haydock method (c) the KPM with Jackson
           damping (d) the KPM without Jackson damping for $M=100$.  }
  \label{fig:dos_lap2d}
\end{figure}

In Fig.~\ref{fig:lap2d_all_blur}, we compare all methods presented in the 
previous section.  We observe that the Lanczos method seems to outperform all 
other methods, especially when $M$ is relatively small.  
The use of Jackson damping in KPM does not
appear to improve the accuracy of the approximation. To some extent, this
is not surprising because the true DOS has many sharp peaks (see
Fig.~\ref{fig:dos_lap2d}) even after 
it is regularized.  Hence, using Jackson damping, which tends to
over-regularize the KPM approximation, may not be able to capture these
sharp peaks.  The DGL method, the KPM and KPML behave similarly, as is 
expected. 

In Fig.~\ref{fig:dos_lap2d}, we compare $\phi_{\sigma}(t)$ and $\tphi(t)$ 
directly for the Lanczos, Haydock, KPM with and without Jackson damping.  
To see the accuracy of different methods more clearly, we 
choose a higher resolution by setting $\sigma$ and $\eta$ to $0.05$.
We note that the meaning of $\tphi(t)$ is different for different
methods.  For Lanczos, $\tphi(t)$ is the approximate DOS obtained using
the Gaussian blurring.  For Haydock, $\tphi(t)$ is the approximate DOS obtained using
the Lorentzian Gaussian blurring. For KPM (with and without Jackson
damping), we first evaluate $\tphi(t)$ as in section~\ref{sec:kpm}, and
then plot instead the quantity $\average{\wt{\phi}(\cdot),
\gsz(\cdot-t)}$.  In this sense, the exact and approximate DOS are
regularized on the same footing. The same procedure is adopted for other
numerical examples in this section as well.

We use $M=100, \nvec=100$ for all methods.  In this case, a visual 
inspection of the approximate DOS plots in Fig.~\ref{fig:dos_lap2d} yields 
the same conclusion that we reached earlier based on the measured errors 
shown in Fig.~\ref{fig:lap2d_all_blur}.  Lanczos appears to be the most 
accurate among all methods.  The DOS curves generated from both the Lanczos and
the Haydock methods are above zero.  The peaks in the DOS curve produced
by the Haydock method are not as sharp as those produced by the Lanczos
method.  This is because Haydock uses a Lorentzian to regularize the
Dirac $\delta$-function, whereas the Lanczos method uses a Gaussian
function to blur the Dirac $\delta$-function centered at Ritz values.  The KPM
method without Jackson damping does not preserve the non-negativity of the
approximate DOS, and Gibbs oscillation is clearly observed in
Fig.~\ref{fig:dos_lap2d} (d).  Finally, KPM with Jackson damping
preserves the non-negativity of the approximate DOS.  However, the use of
Jackson damping leads to missing several peaks in the DOS, as is
illustrated in Fig.~\ref{fig:lap2d_all_blur}. The behavior of DGL
and KPML are similar to that of
KPM without Jackson damping.

\subsection{Other test matrices}\label{sec:othermatrix}

In this section, we compare different DOS approximation methods
for two other matrices taken from the Univerisity of Florida Sparse Matrix
collection \cite{FloridaMatrix}. The {\em pe3k} matrix originates from 
vibrational mode calculation of a polyethylene molecule with 3,000 atoms 
\cite{YangNCA2001}. The {\em shwater} matrix originates from a computational 
fluid dynamics simulation.  The size of the {\em pe3k} matrix is 9,000, and the 
size of shwater matrix is 81,920. These two test matrices have quite 
different characteristics in their DOS. The spectrum of the {\em pe3k} matrix
contains a large gap as well as many peaks. The DOS of the {\em shwater} 
matrix is relatively smooth as we will see below.

We set $\sigma$ to $0.3$ in tests presented in this section.
We observe that KPM with Jackson damping only becomes
accurate when the degree of the expanding polynomials ($M$) is high enough, 
and the convergence with respect to $M$ is rather slow.  The DGL method, KPM 
without Jackson damping and KPML behave similarly.  


For the {\em shwater} matrix, which has a relatively smooth spectral density,
the Lanczos method is still the most accurate. It only take $M=50$
Lanczos steps to reach $10^{-3}$ accuracy.  The KPM (without
Jackson damping) and KPML, as well as the DGL method all
require $M>110$ terms to reach the same level of accuracy.  

\begin{figure}[h]
  \begin{center}
    \includegraphics[width=0.4\textwidth]{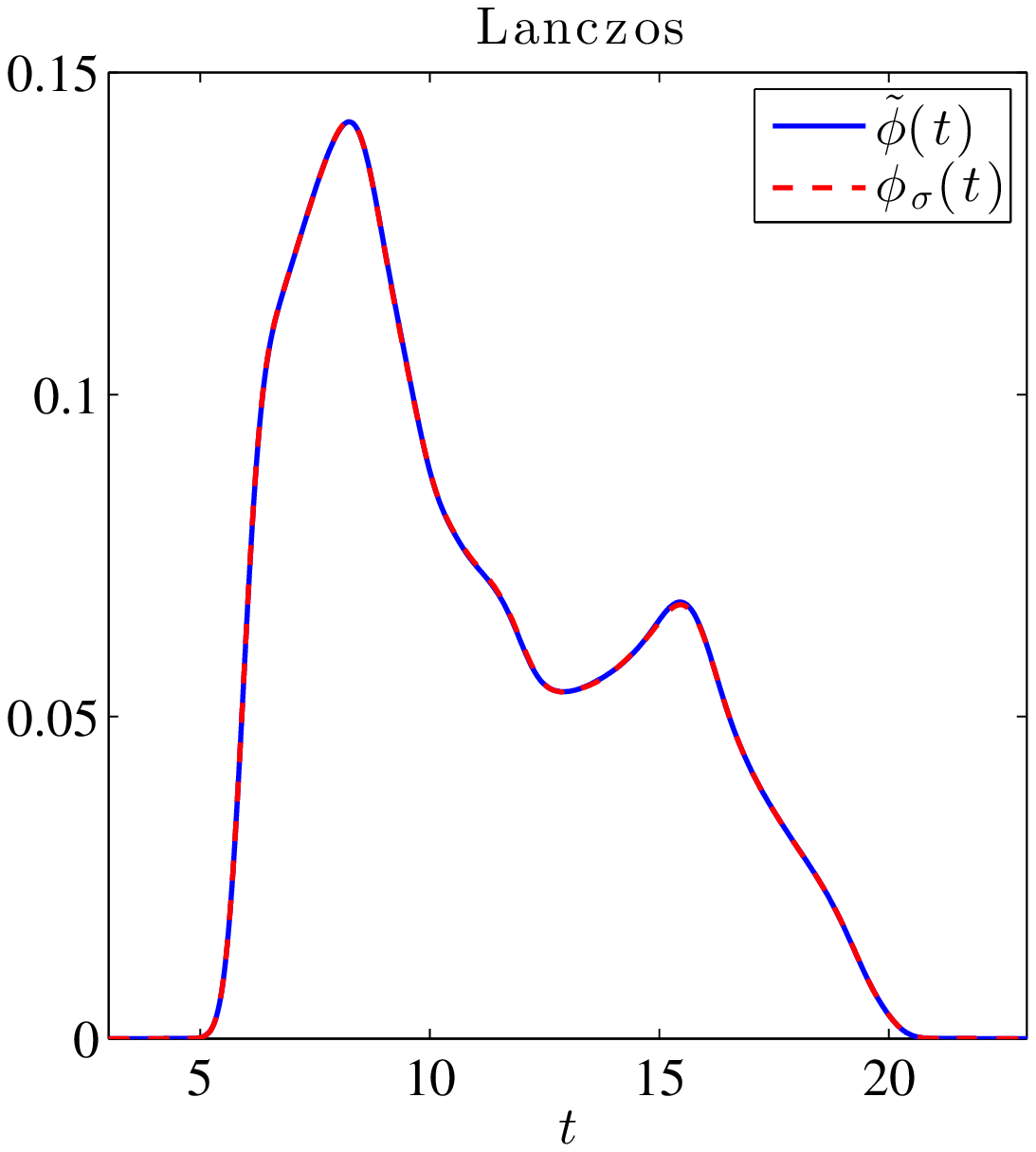}
    \includegraphics[width=0.4\textwidth]{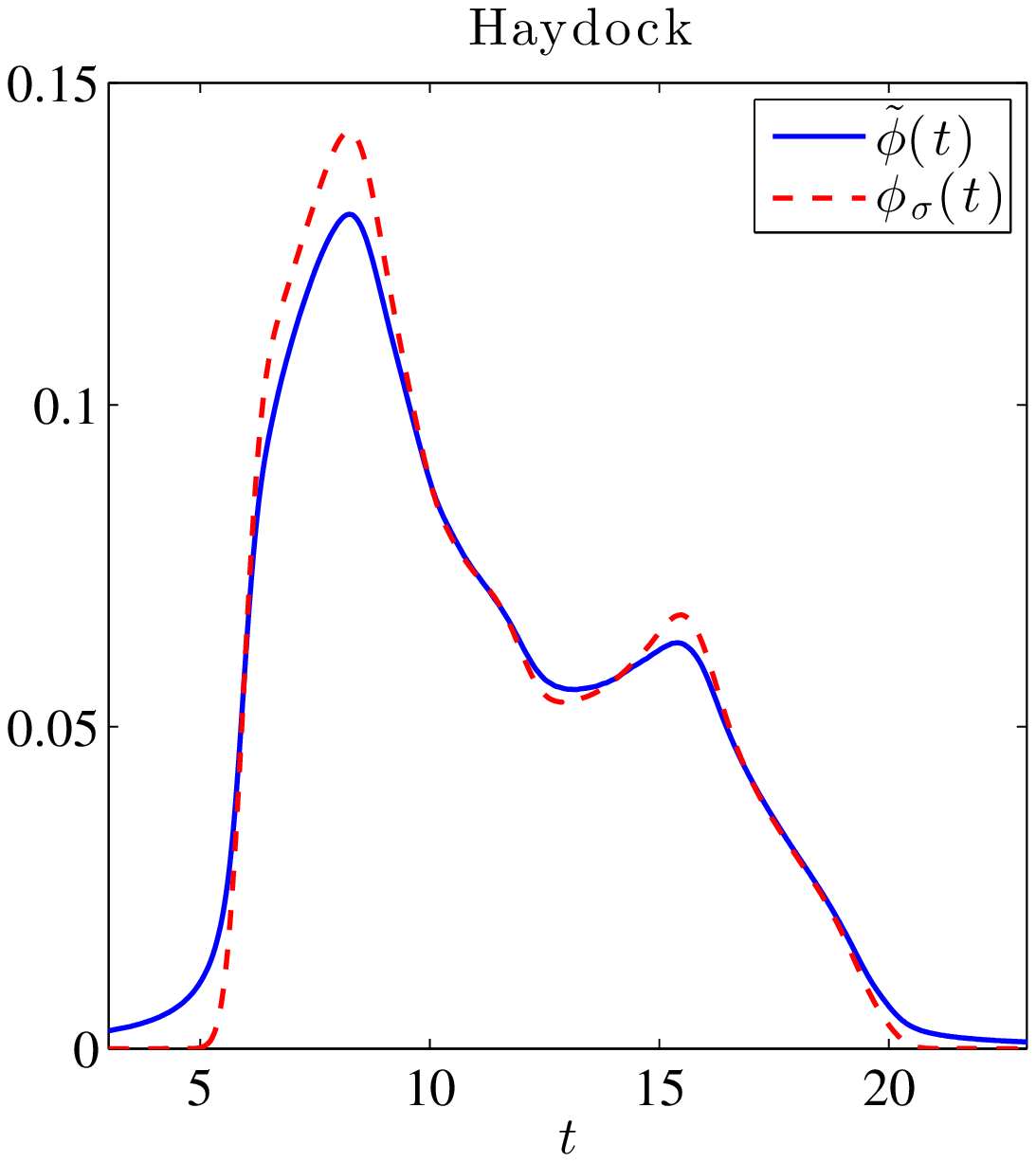} 
    
    \includegraphics[width=0.4\textwidth]{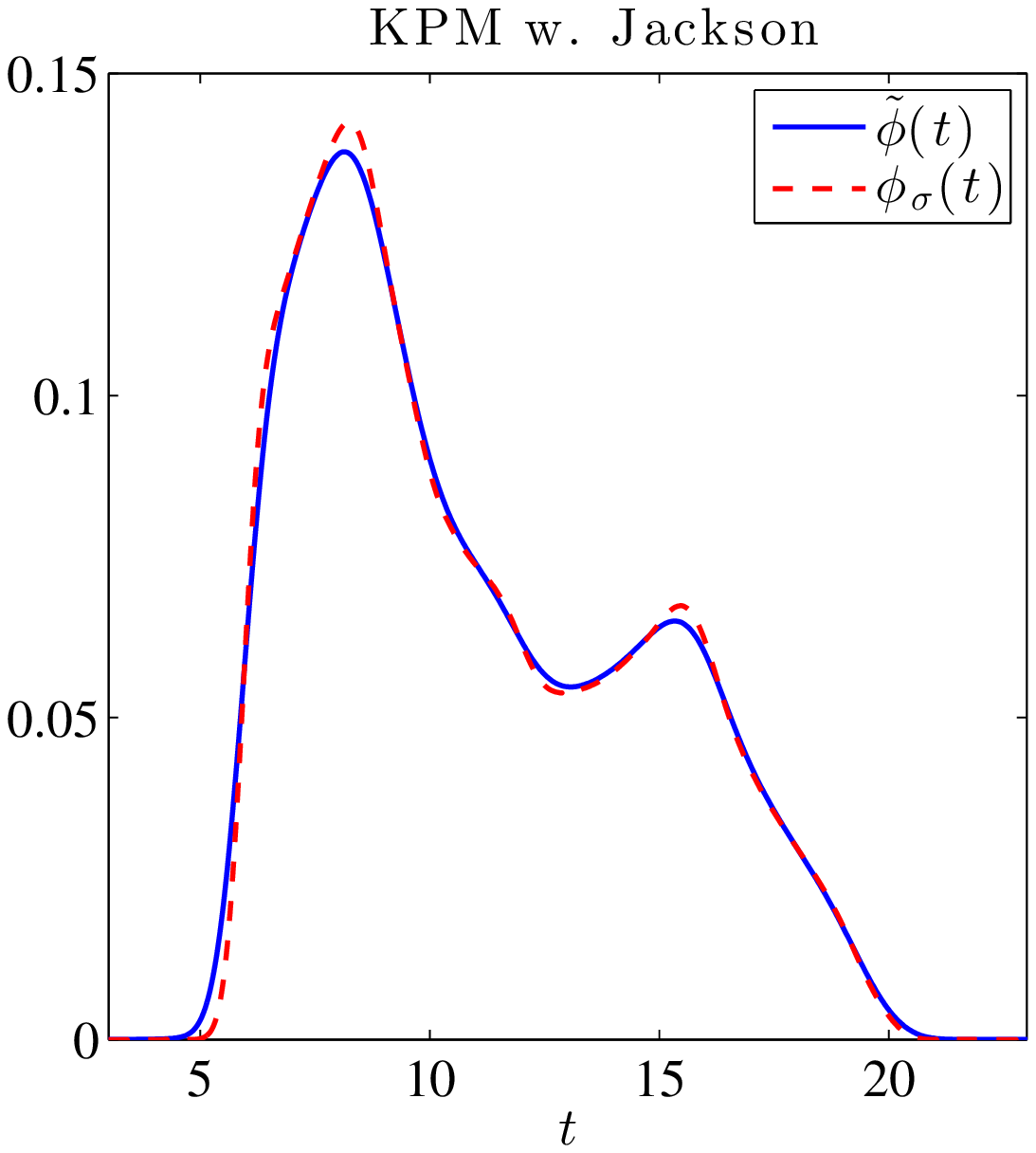}
    \includegraphics[width=0.4\textwidth]{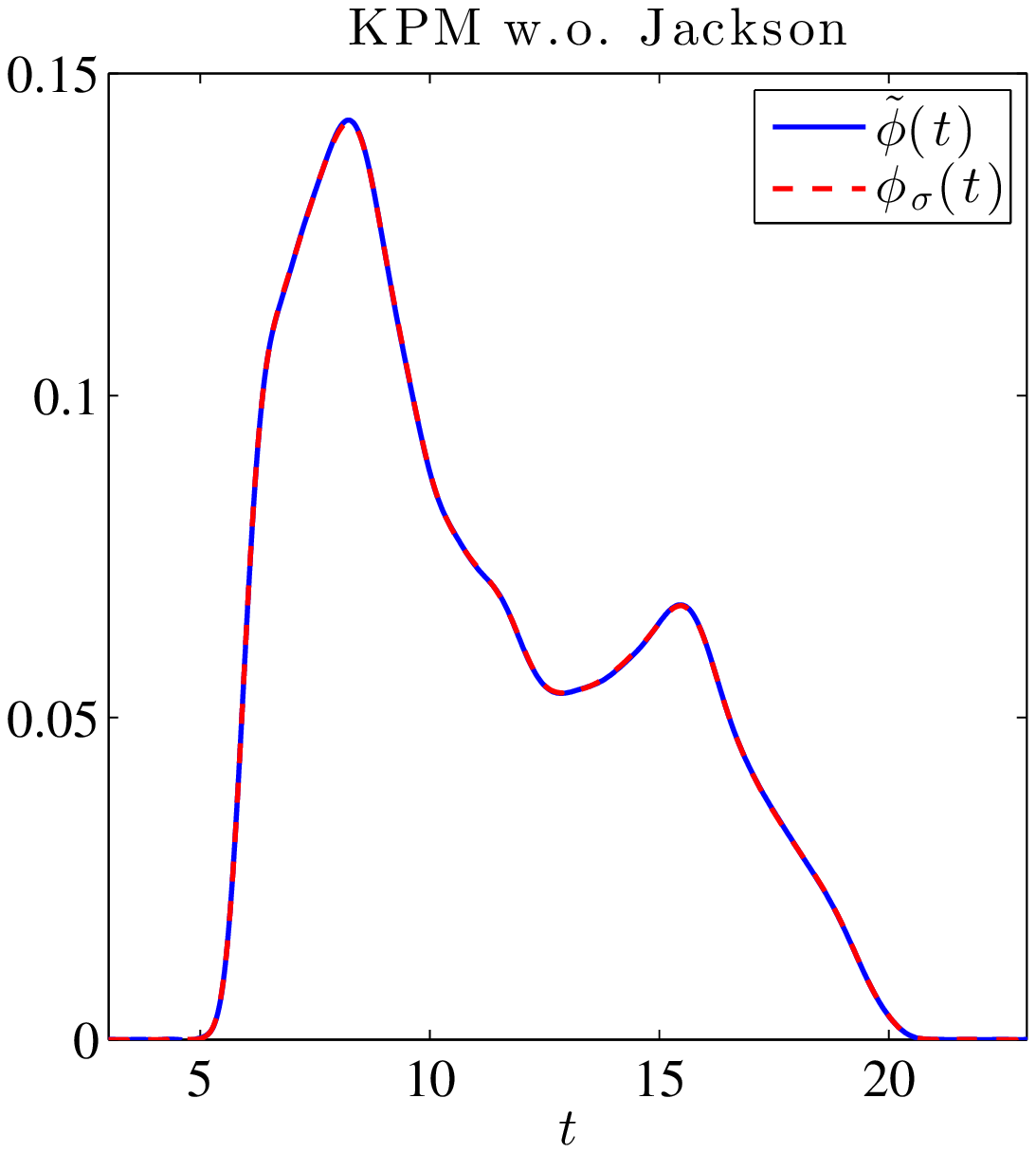}
  \end{center}
  \caption{Comparing the regularized DOS (with $\sigma =0.3$) with
  the approximate DOS produced by (a) the Lanczos method 
  (b) the Haydock method (c) the KPM with Jackson damping 
  (d) the KPM without Jackson damping for the {\em shwater} matrix.}
  \label{fig:dos_shwater}
\end{figure}

Fig.~\ref{fig:dos_shwater} shows that when we set $M=100$, $\nvec=100$, 
the KPM without Jackson damping yields accurate approximation to the DOS, 
whereas the Jackson damping introduces slightly larger error near the locations
of the peaks and valleys of the DOS curve. This error is due to the use of 
extra smoothing.  The DOS generated by the Haydock method also has larger 
error near the peaks and valleys of the DOS curves. This is due to 
the use of Lorentzian regularization.  

\begin{figure}[h]
  \begin{center}
    \includegraphics[width=0.4\textwidth]{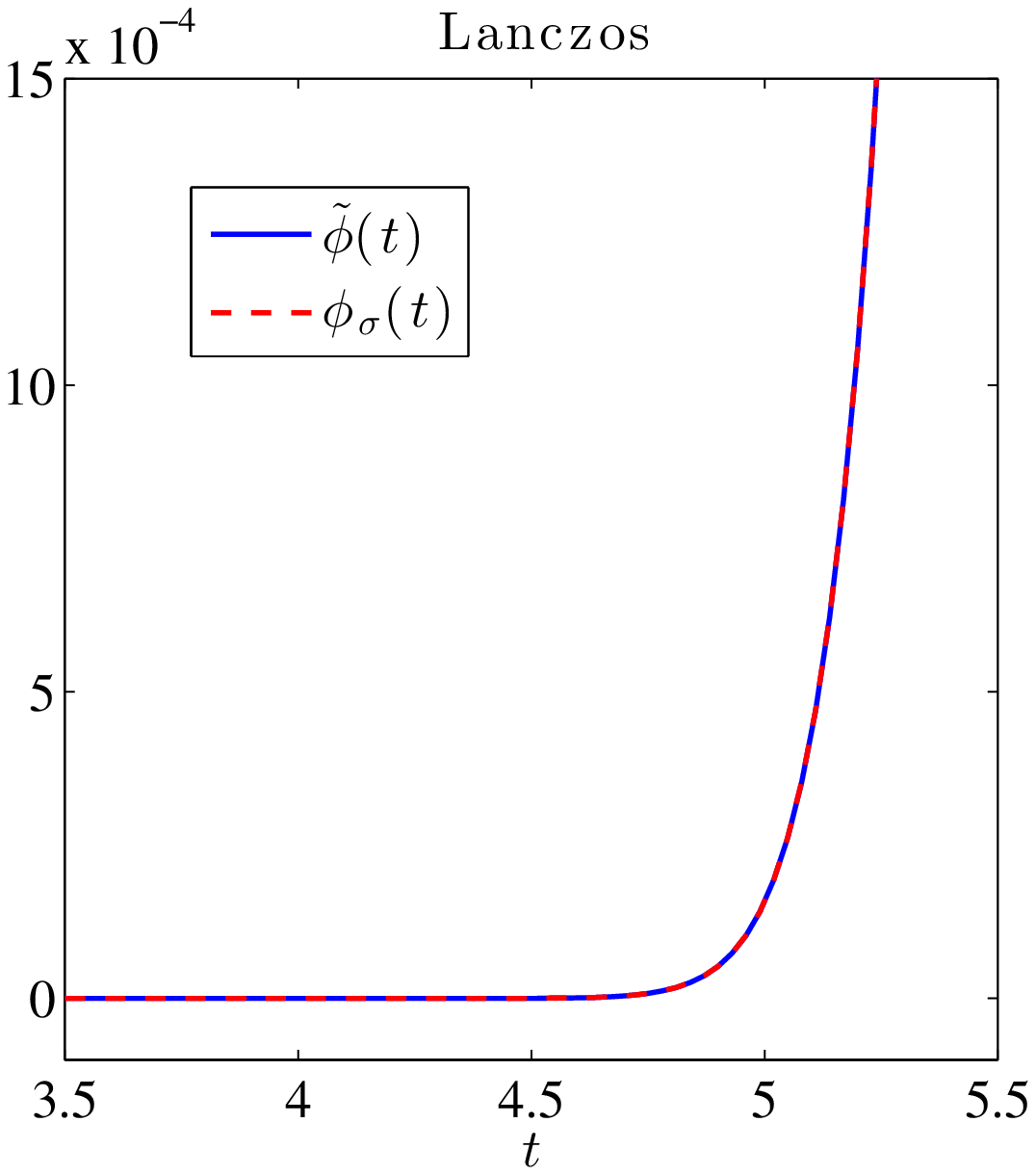}
    \includegraphics[width=0.4\textwidth]{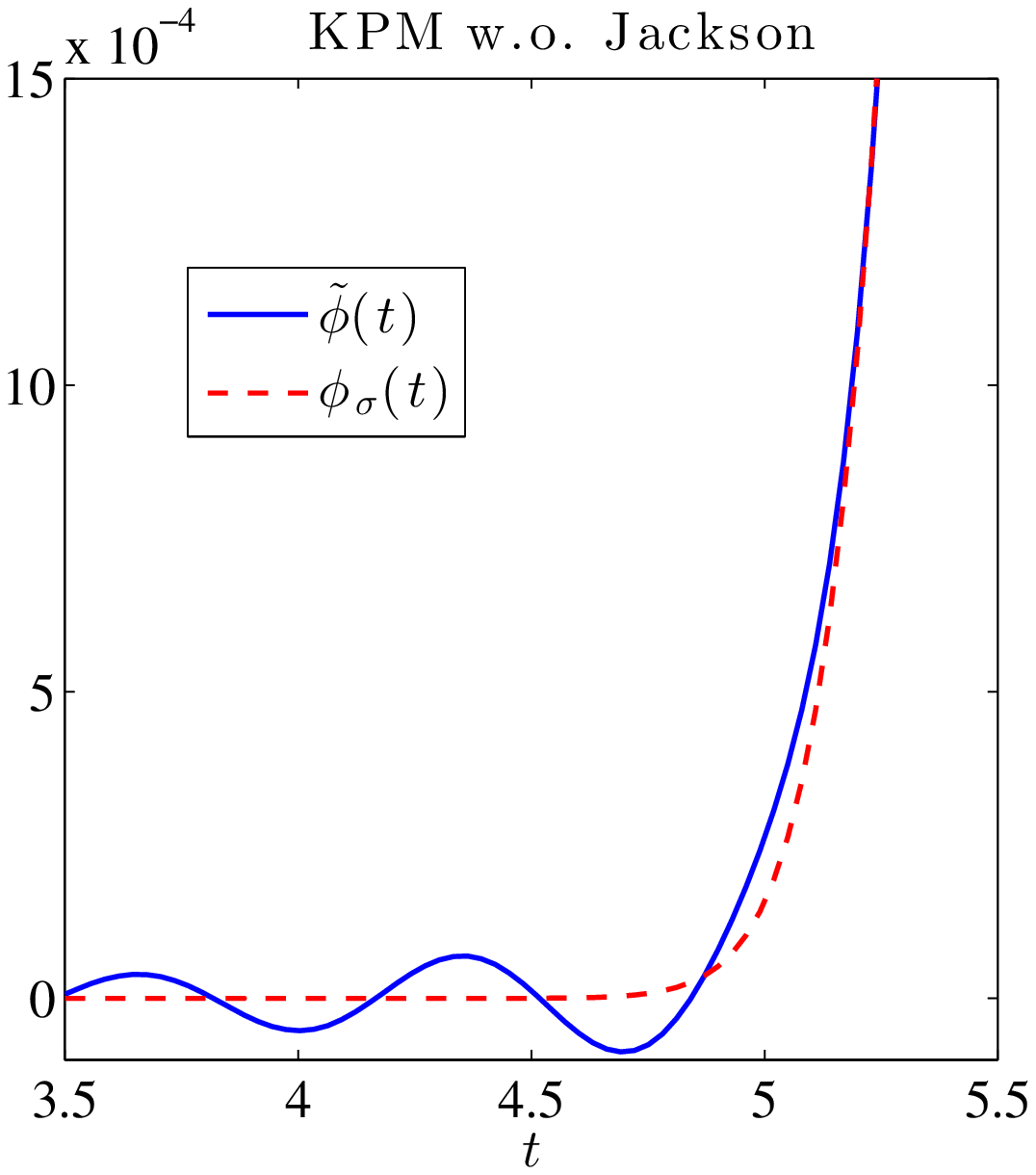}
  \end{center}
  \caption{A comparison of approximation errors of the Lanczos method (a)
  with that of the KPM without Jackson damping (b) at the higher end 
  of the spectrum of the {\em shwater} matrix.}
  \label{fig:dos_shwater_zoom}
\end{figure}

In Fig.~\ref{fig:dos_shwater_zoom},
we zoom into the tail of the DOS curves produced by the Lanczos
and KPM.  It can be seen that Lanczos preserves the non-negativity of the
DOS, whereas KPM does not.  However, since the DOS is smooth, the Gibbs
oscillation is very small, and can only be seen clearly at the tail of the
DOS curve.


\begin{figure}[h] \begin{center}
    \includegraphics[width=0.4\textwidth]{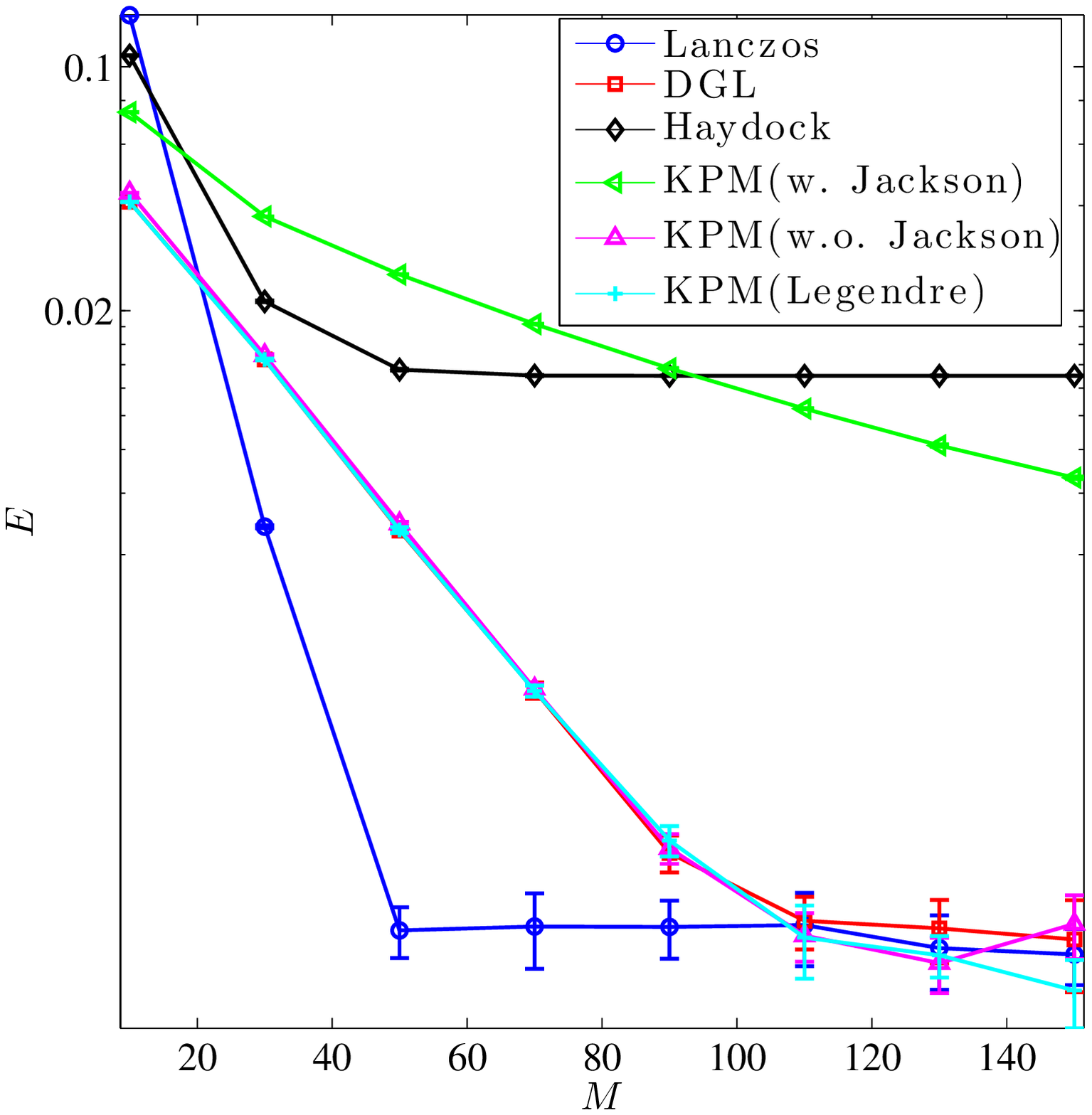}
  \end{center} 
  \caption{A comparision of approximation errors of all DOS approximation
           methods applied to the {\em shwater} matrix for different $M$ values.
           The regularization parameter $\sigma$ is set to 0.3.}
  \label{fig:shwater_all_blur}
\end{figure}

\begin{figure}[h]
  \begin{center}
    \includegraphics[width=0.4\textwidth]{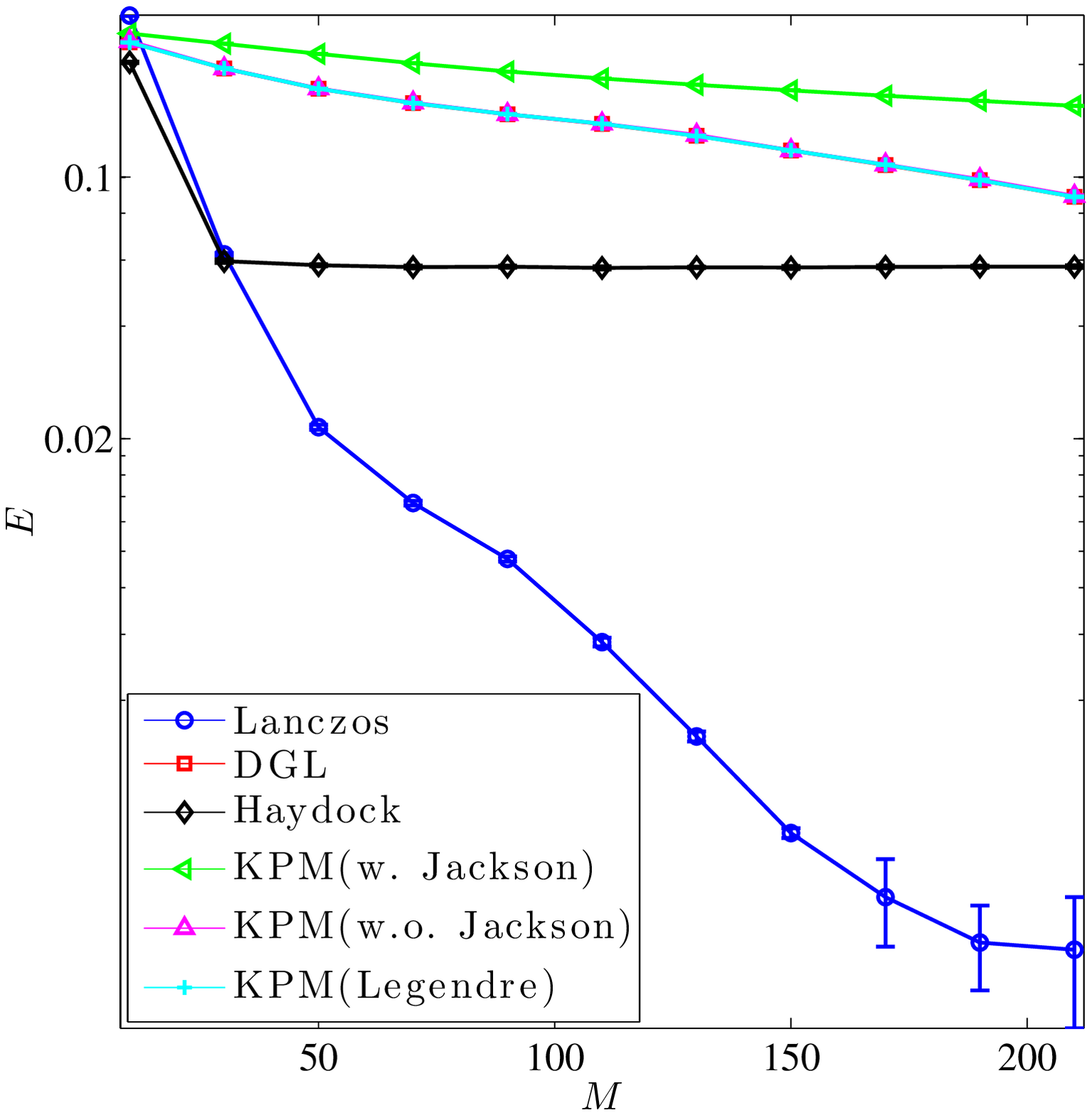}
  \end{center} 
  \caption{A comparison of approximation errors of all DOS approximation
           methods applied to the {\em pe3k} matrix for different $M$ values.
           The regularization parameter $\sigma$ is set to 0.3.}
  \label{fig:pe3k_all_blur}
\end{figure}

\begin{figure}[h]
  \begin{center}
    \includegraphics[width=0.4\textwidth]{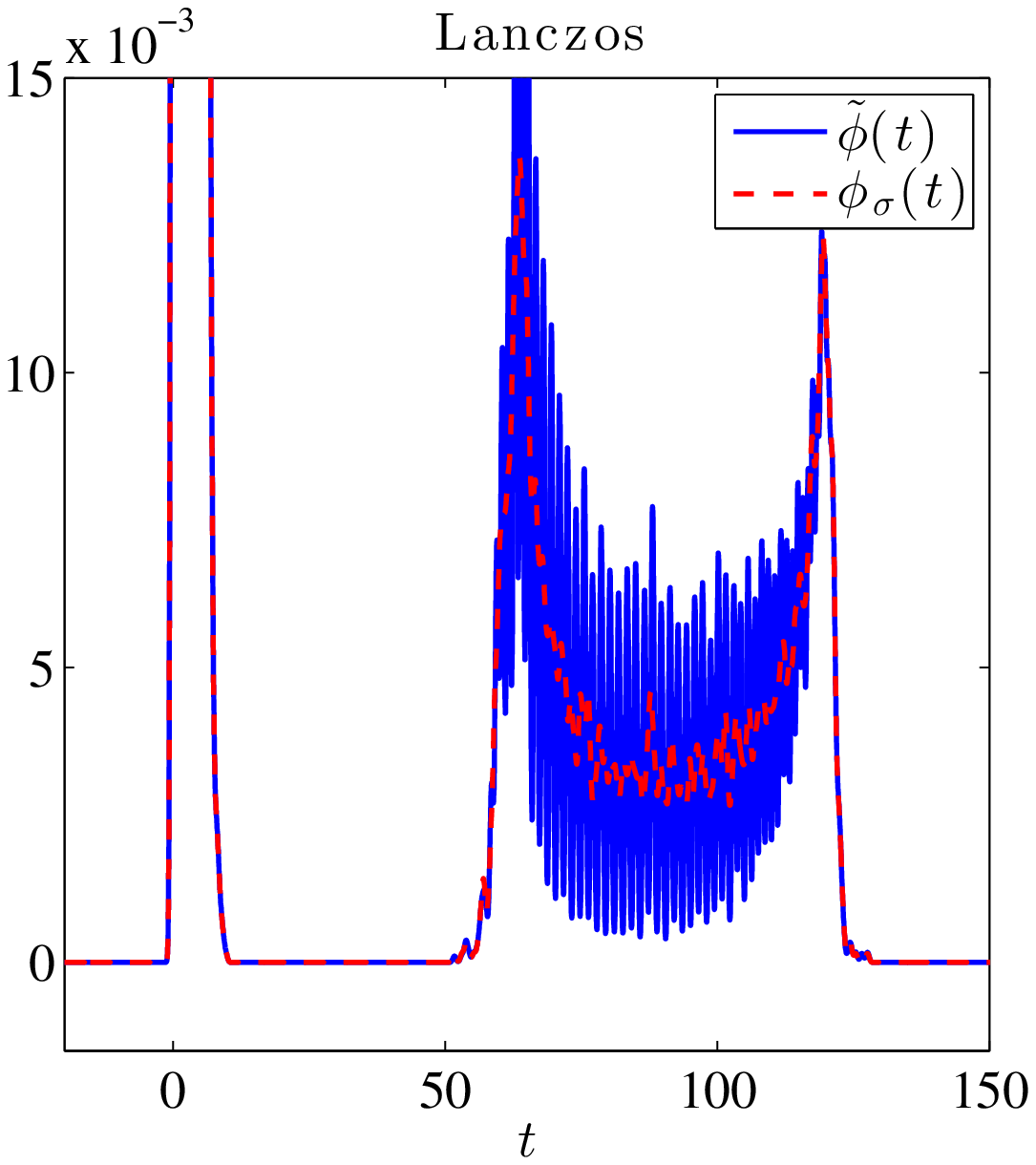}
    \includegraphics[width=0.4\textwidth]{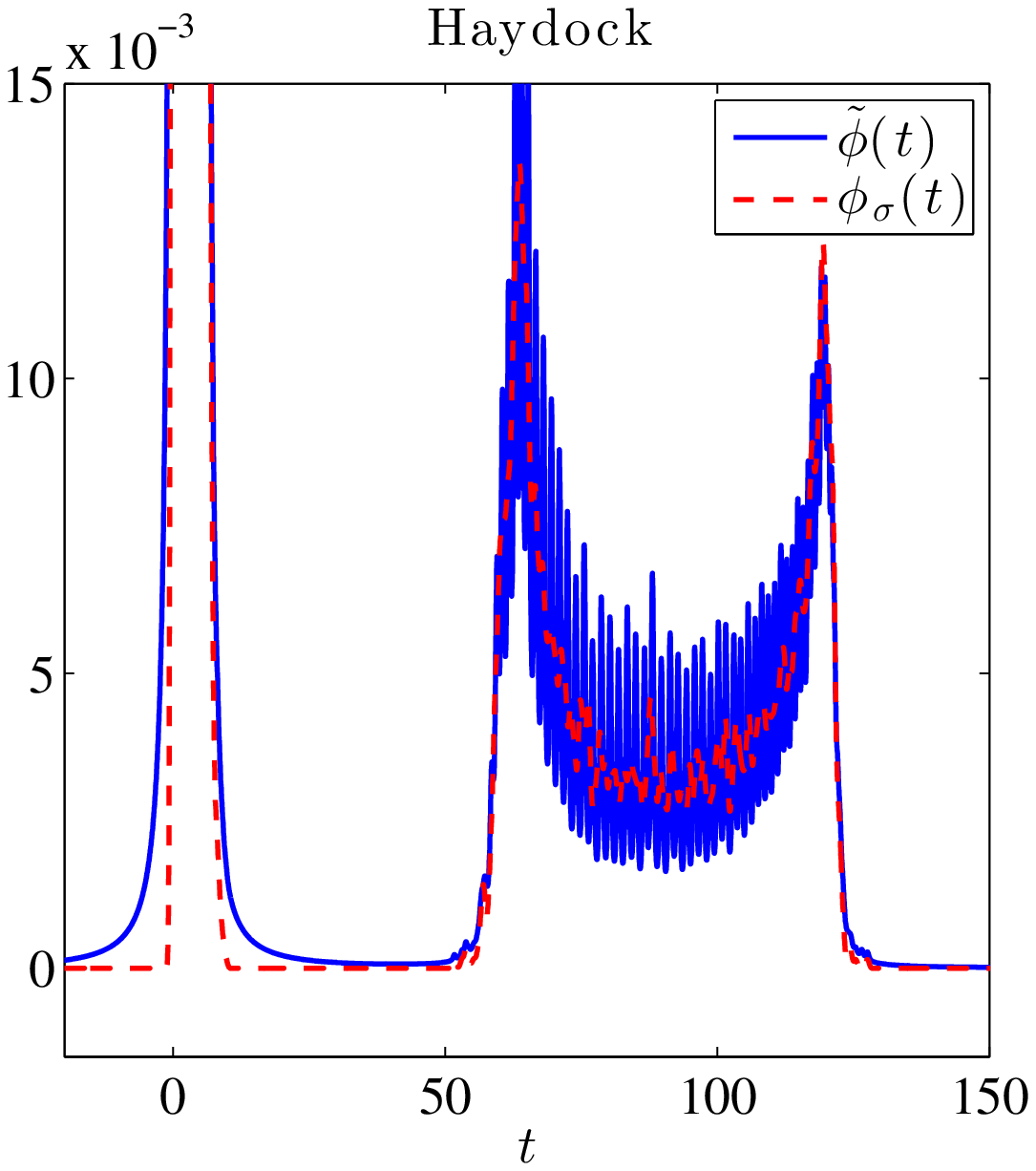} 
    
    \includegraphics[width=0.4\textwidth]{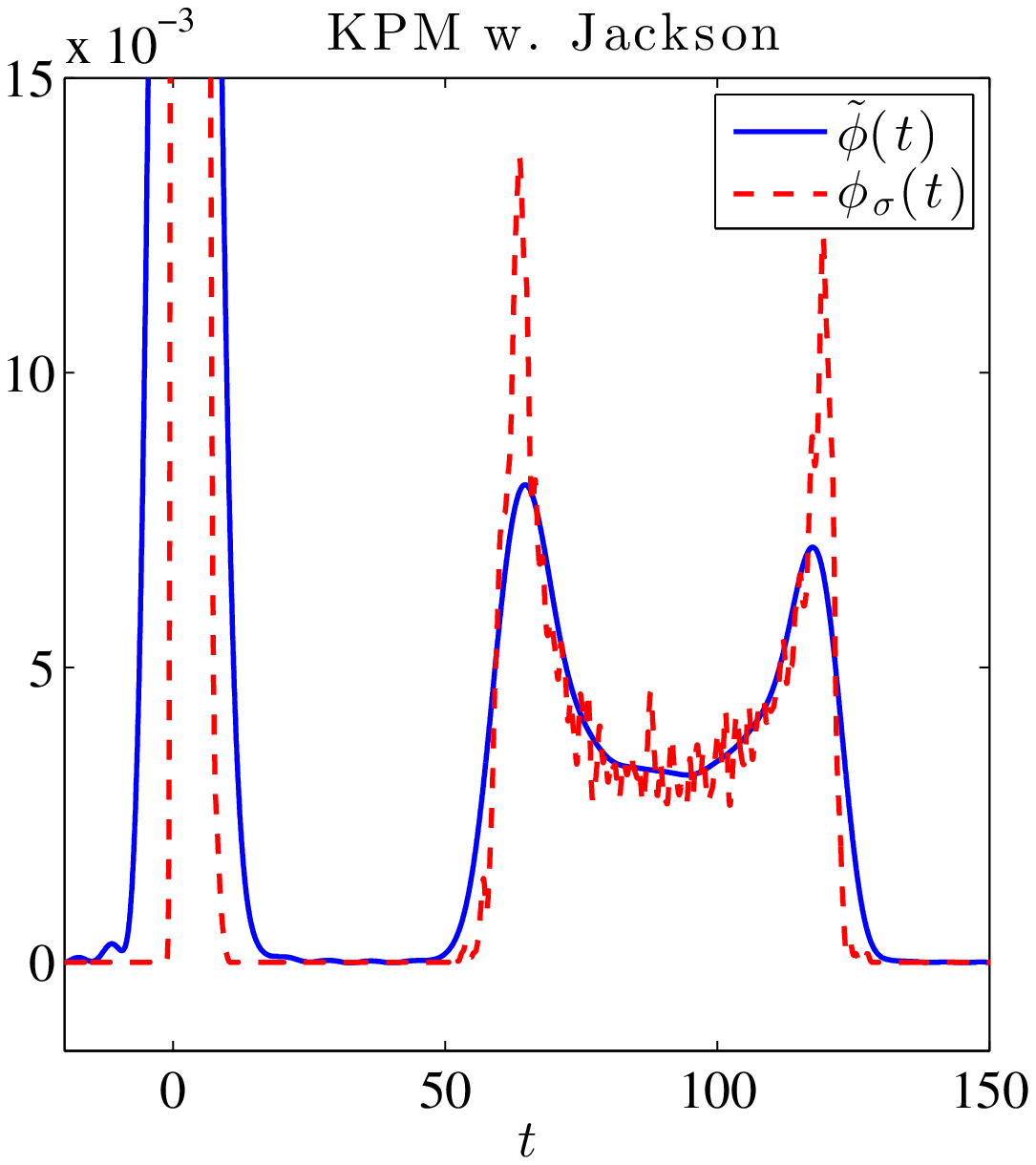}
    \includegraphics[width=0.4\textwidth]{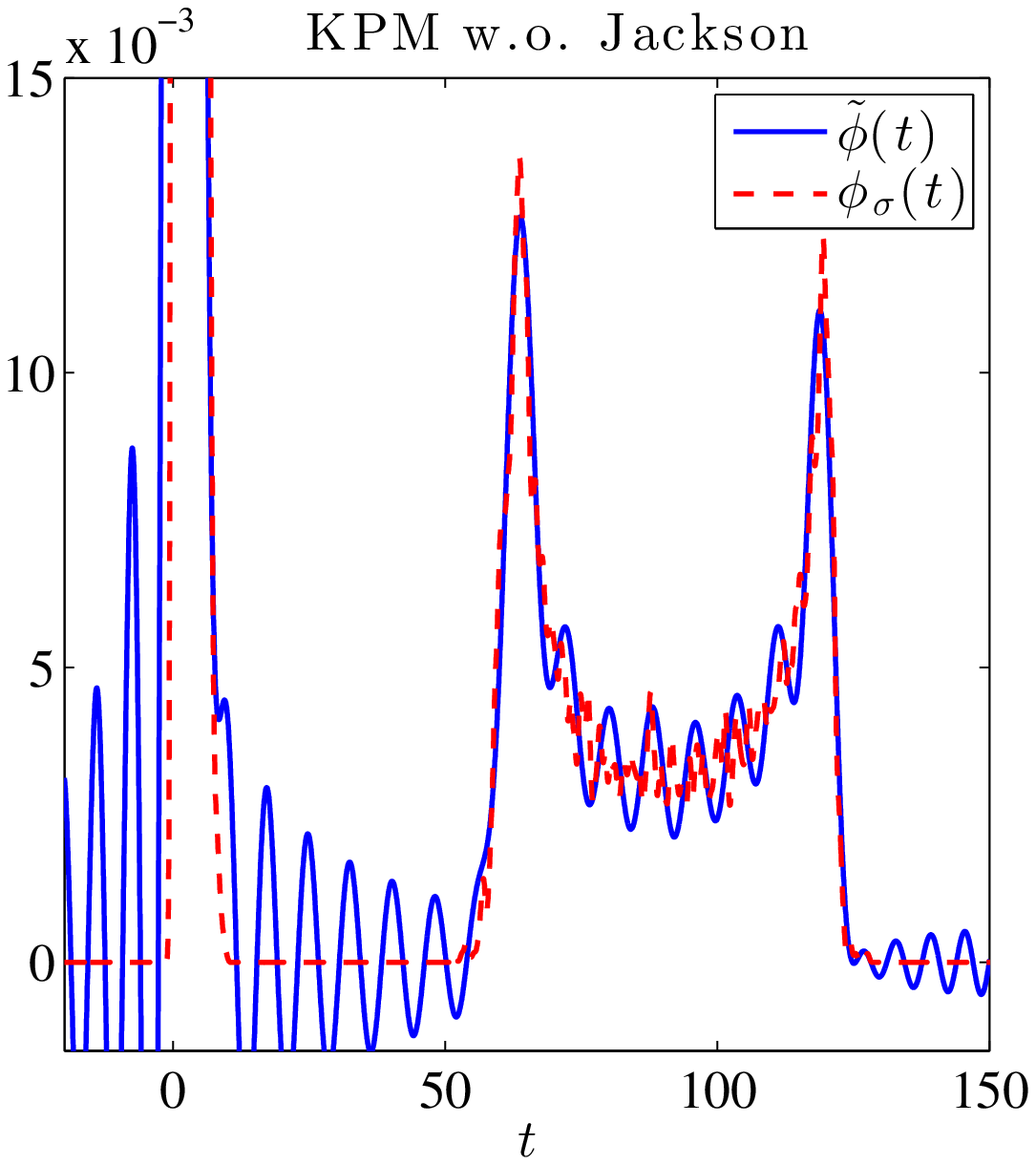}
  \end{center}
  \caption{Comparing the regularized DOS (with $\sigma =0.3$) with
  the approximate DOS produced by (a) the Lanczos method 
  (b) the Haydock method (c) the KPM with Jackson damping 
  (d) the KPM without Jackson damping for the {\em pe3k} matrix.}
  \label{fig:dos_pe3k}
\end{figure}

For the {\em pe3k} matrix, 
Fig.~\ref{fig:pe3k_all_blur} shows that the Lanczos method 
is significantly more accurate than other methods, followed by 
the Haydock method.  This difference in accuracy can be further observed 
in Fig.~\ref{fig:dos_pe3k}, which compares the regularized DOS with
the approximate DOS for the Lanczos method, the Haydock method, and the 
KPM with and without Jackson damping.  We use $M=100$ and $\nvec=100$. 
The {\em pe3k} matrix has a large gap between the low and high ends of the 
spectrum. Without the use of Jackson damping, the KPM produces large 
oscillations over the entire spectrum.  We observed similar
behavior for DGL and KPML.
Adding Jackson damping reduces oscillation in the approximate DOS. 
However, it leads to an over-regularized DOS approximation and is not
accurate.  

\subsection{Application: Heat capacity calculation}\label{sec:application}

At the end of section~\ref{sec:err1} we discussed that there are
different ways for regularizing the DOS depending on the applications.
Here we give an example of the calculation of the heat capacity of a
molecule.  The heat capacity is a thermodynamic property and is
defined as~\cite{mcquarrie,heatcap}
\begin{equation}
C_v =  \int_0^\infty k_B\frac{(\hbar \omega c/k_BT)^2 e^{-\hbar \omega
c/k_B T} }{(1-e^{-\hbar \omega c/k_B T})^2}  \phi(\omega) d\omega,
\label{eq:heatcv}
\end{equation}
where $k_B$ is the Boltzmann constant, $c$ is the speed of light,
$\hbar$ is Planck's constant, $T$ is the temperature and $\omega =
\sqrt{\lambda}$ is the vibration frequency.  

Here, if we define
\begin{equation}
  g(\omega) = k_B\frac{(\hbar \omega c/k_BT)^2 e^{-\hbar \omega
  c/k_B T} }{(1-e^{-\hbar \omega c/k_B T})^2},
  \label{eqn:gfunc}
\end{equation}
and define the DOS $\phi(\omega)$ using the square root of the
eigenvalues of the Hessian associated with a molecular potential
function with respect to atomic coordinates of the molecule, we have
\[
C_{v} = \average{\phi,g}.
\]
Therefore the error can be measured directly using
Eq.~\eqref{eqn:errormetric_g}.

In the following, we take the Hessian to be the modified
Laplacian matrix and the {\em pe3k} matrix, and compute the 
corresponding heat capacity $C_v(T)$ for different temperature values $T$.
We note that here the computed values of $C_v(T)$ do not carry any physical 
meaning.  They merely serve as a proof of principle for assessing the 
accuracy of the estimated DOS.  We compare the KPM and the Lanczos method.
All computations are done using $M=40$ MATVECs. Each computed $C_v$ is
an averaged value over $100$ runs. To facilitate the comparison, we
normalize $C_v$ so that its maximum value is $1$.  The Lanczos
method is also fully flexible when the error metric is changed. To this
end we regularize the distribution obtained from Ritz values not by
Gaussians, but by the function $g$ in this application. In other words, in
Eq.~\eqref{eqn:landos} we replace $g_{\sigma}$ by the function $g$ in
Eq.~\eqref{eqn:gfunc}.

\begin{figure}[htb]
	\begin{center}
	\includegraphics[width=0.40\textwidth]{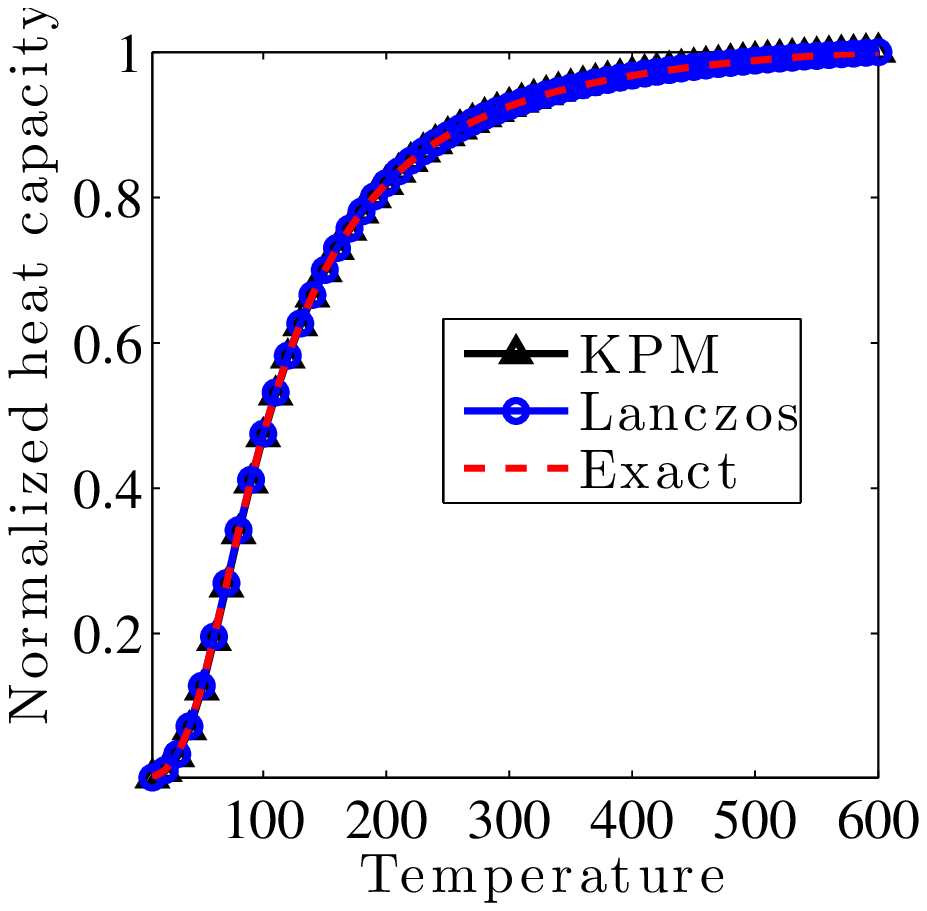}\quad
	\includegraphics[width=0.42\textwidth]{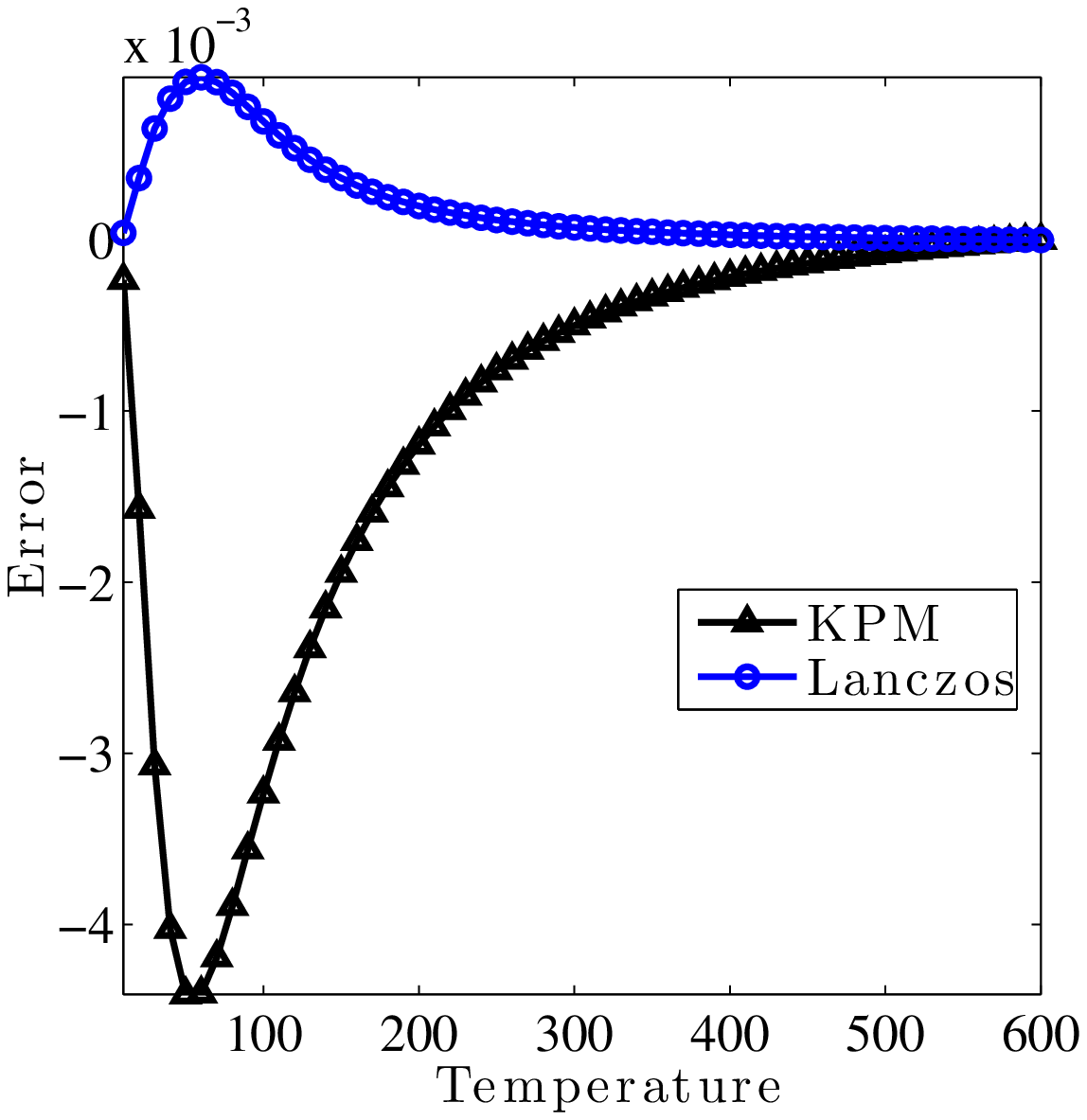}\quad
        \end{center}
  \caption{A comparison of the approximate heat capacity produced by the 
           Lanczos and the KPM with the ``exact" heat capacity at 
           different temperatures (left) , and the approximation errors 
           produced by these methods at different temperature values (right)
           for the modified Laplacian.}
	\label{fig:heatlap2d}
\end{figure}

Fig.~\ref{fig:heatlap2d} shows that both the KPM and the Lanczos method 
correctly reproduce the normalized $C_v(T)$ for the modified Laplacian
matrix. We also plot the error generated in both the KPM and
the Lanczos method.  We observe that the error associated with the 
Lanczos method is slightly smaller.  This observation agrees with previous 
results that demonstrate the effectiveness and accuracy of both the 
KPM and the Lanczos methods for computing a relatively smooth DOS.

\begin{figure}[htb]
	\begin{center}
	\includegraphics[width=0.41\textwidth]{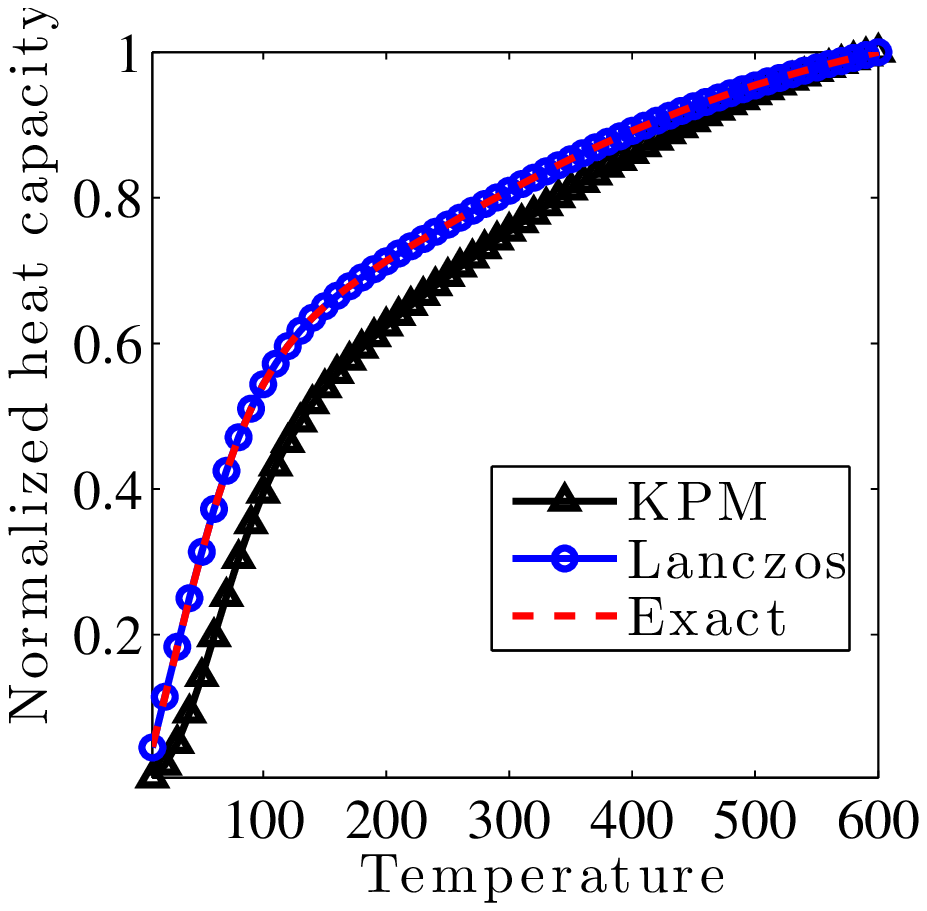}\quad
	\includegraphics[width=0.40\textwidth]{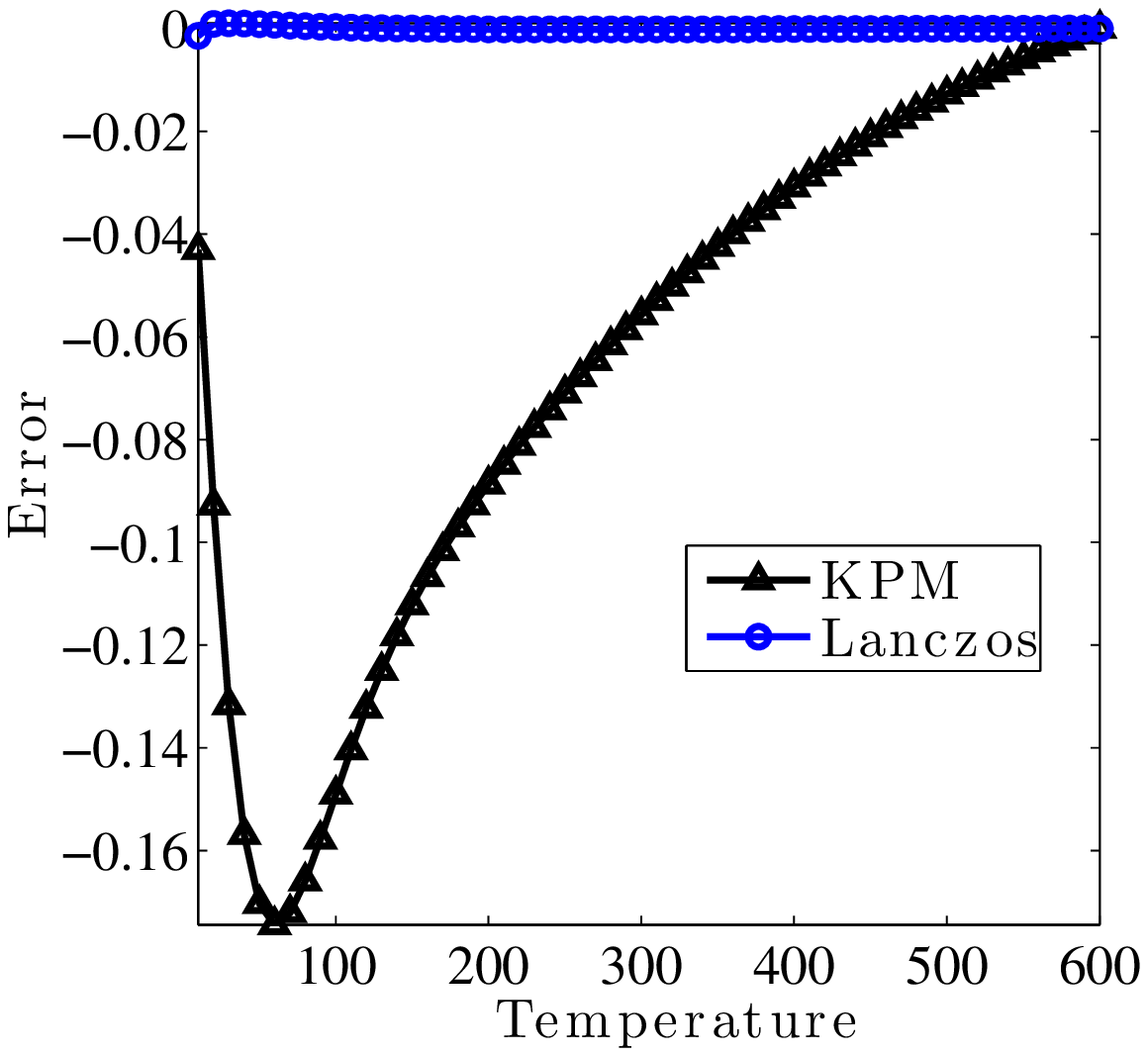}\quad
        \end{center}
  \caption{A comparison of the approximate heat capacity produced by the 
           Lanczos and the KPM with the ``exact" heat capacity at 
           different temperatures (left) , and the approximation errors 
           produced by these methods at different temperature values (right)
           for the {\em pe3k} matrix.}
	\label{fig:heatpe3k}
\end{figure}

Fig.~\ref{fig:heatpe3k} shows that, for the {\em pe3k} matrix, the KPM 
approximation of $C_v(T)$ exhibits much larger error than that produced 
by the Lanczos method.  This observation agrees with the results 
shown in Figure~\ref{fig:dos_pe3k}, which suggests that the Lanczos method
yields a much more accurate DOS estimate, especially when $M$ is relatively
small.


\section{Conclusion}
We surveyed numerical algorithms for estimating the spectral
density of a real symmetric matrix $A$ from a numerical linear algebra
perspective. The algorithms can be categorized into two classes.
The first class contains the KPM method an its variants.
The KPM is based on constructing polynomial approximations to 
Dirac $\delta$-``functions" or regularized $\delta$-``functions".  
We showed that the Lanczos spectroscopic
method is equivalent to KPM even though it is derived from different
view points.  The DGL method is slightly different, but can be
viewed as a polynomial expansion of a regularized spectral density. 
It is more flexible because it allows
polynomials of different degrees to be used at different spectral
locations.

The second class of methods is based on the classical Lanczos
procedure for partially tridiagonalizing $A$.  Both the Lanczos 
and the Haydock methods make use of eigenvalues and eigenvectors
of the tridiagonal matrix to construct approximations to the DOS. 
They only differ in the type of regularization they use to interpolate 
spectral density from Ritz values to other locations in the spectrum.
The Lanczos method uses a Gaussian blurring function, whereas the
Haydock method uses a Lorentzian. Because a Lorentzian decreases to
zero at a much slower rate than a Gaussian away from its peak, 
it is less effective when a high resolution spectral density is 
needed.

Regularization through the use of Gaussian blurring of 
$\delta$-``functions" not only allows us to specify the desired resolution
of the approximation, but also allows us to properly define
an error metric for measuring the accuracy of the approximation 
in a rigorous and quantitative manner.

The KPM and its variants require estimating the trace of 
$A$ or $p(A)$ where $p(t)$ is a polynomial.  
An important technique for obtaining such an estimate is the stochastic 
sampling and averaging of the Rayleigh quotient $v_0^T p(A) v_0 /
v_{0}^{T} v_{0}$.  Averaging the tridiagonal matrices
produced by the Lanczos procedure started from randomly generated 
starting vectors ensures that the approximation contains
equal contributions from all spectral components of $A$.
This is an important requirement of the Lanczos and Haydock algorithms.

Our numerical tests show that the Lanczos method consistently outperforms
other methods in term of the accuracy of the approximation, especially
when a few MATVECs are used in the computation.  
Furthermore, both the Lanczos and Haydock
algorithms guarantee that the approximate DOS is non-negative. This is 
a desirable feature of any DOS approximation.  Another nice property
of the Lanczos and Haydock algorithms is that in the limit of $M=n$, they
fully recovered a regularized DOS.

The KPM and its variants appear to work well when the DOS to be approximated is relatively smooth. 
They are less effective when the DOS contains many peaks or when the
spectrum of $A$ contains large gaps. We found the use of 
Jackson damping can remove the Gibbs oscillation of KPM. However, it
tends to over-regularized the approximate DOS and misses important features (peaks) of the DOS.  

\section*{Acknowledgments}

This work is partially supported by the Laboratory Directed Research and
Development Program of Lawrence Berkeley National Laboratory under the
U.S. Department of Energy contract number DE-AC02-05CH11231 (L. L. and
C. Y.), and by Scientific Discovery through Advanced Computing (SciDAC)
program funded by U.S.  Department of Energy, Office of Science,
Advanced Scientific Computing Research and Basic Energy Sciences
DE-SC0008877 (Y. S. and C. Y.).

\FloatBarrier

\appendix

\section{Further discussion on the KPM method}\label{sec:app_kpm}

For KPM, a common approach used to damp the Gibbs
oscillations is to use the Chebyshev-Jackson 
approximation~\cite{Jackson1930,Rivlin2003,JayKimSaadEtAl1999}, which
modulates the coefficients $\mu_{k}$ with a damping factor $g_{k}^{M}$ defined
by
\begin{equation}
	g_{k}^{M} = \frac{\left(1-\frac{k}{M+2}\right) \sin(\alpha_{M})
	\cos(k\alpha_{M}) + \frac{1}{M+2}\cos(\alpha_{M})\sin(k
	\alpha_{M})}{\sin(\alpha_{M})},
	\label{eq:jackson}
\end{equation}
where $\alpha_{M}=\frac{\pi}{M+2}$. 
Consequently, the damped Chebyshev expansion has the form
\[
\wt{\phi}_{M}(t)  = \sum_{k=0}^M  \mu_k g_{k}^{M} T_k(t).
\]
The approximation of Jackson damping is demonstrated in
Fig.~\ref{fig:delCheb40}.

Another variant can be derived from that
Chebyshev polynomials are not  the only types of orthogonal polynomials
that can be used in the expansion. We can use any other type of orthogonal polynomials.
The only practical requirement is that we explicitly know the 3-term
recurrence for the polynomials. 
For example, we can use the Legendre polynomials $L_k(t)$ which obey
the following 3-term recursion
\[
L_0(t) = 1, \ \ L_1(t) = t, \ \ (k+1)L_{k+1}(t) = (2k+1)t L_k(t) - 
k L_{k-1}(t).
\]
See, for example~\cite{Davis}, for 
three-term recurrences for a wide class of such polynomials,
e.g., all those belonging to the Jacobi class, which include
Legendre and Chebyshev polynomials as particular cases.

From a computational point of view,  some  savings 
in time can be achieved if we are willing to store more vectors.
This is due to the formula:
\[
T_p (t) T_q(t) = \half \left[ T_{p+q} (t) - T_{|p-q|} (t) \right], 
\]
from which we obtain
\[
  T_{p+q} (t) = 
2 \ T_p (t) T_q(t)  + T_{|p-q|} (t) .
\]
For a given $k$ we can use the above formula with 
$p= \lceil k/2 \rceil$ and $q=k-p$.
This requires that we compute and store $v_r = T_r(A) v_{0} $ for
$r \le p$. Then the moments $v^T_{0} T_r(A) v_{0}$ for $r\le p$ can be
computed in the usual way, and for $r=p+q>p$ we can use the formula:
\[
v_{0}^T T_{p+q} (A) v_{0} = 
2 \ v_p ^T v_q + v^T_{0} v_{|p-q|} .
\]
This saves 1/2 of the matrix-vector products at the expense of
storing all the previous $\{v_r\}$, and therefore it is not practical for high
degree  polynomials.

\section{Details on the derivation of the DGL method}\label{sec:app_dgl}

We now calculate the $\gamma_k$'s starting with $\gamma_0$. 
Since $L_0(\lambda)=1$, a change of variable $t \gets (s-t)/\sqrt{2 \sigma^2}$ 
yields
\eq{eq:gam0} \gamma_0 
= 
\sigma\sqrt{\frac{\pi}{2}}
\left[ \mbox{erf} \left(\frac{1-t}{\sqrt{2}\sigma}\right) -
\mbox{erf} 
\left(\frac{-1-t}{\sqrt{2}\sigma}\right) \right]  = 
\sigma\sqrt{\frac{\pi}{2}}
\left[\mbox{erf} \left(\frac{1-t}{\sqrt{2}\sigma} \right)+
\mbox{erf} \left(\frac{1+t}{\sqrt{2}\sigma}
\right) \right]  ,
\en
where we have used the standard error function: 
\[
\mbox{erf} (x) = \frac{2}{\sqrt{\pi}} \int_0^x e^{-t^2} dt \ .
\] 

Now consider a general coefficient $\gamma_{k+1}$ with $k \ge 0$.
There does not seem to exist a closed form formula for $\gamma_k$ for
a general $k$.  
However, these coefficients can be obtained by a recurrence relation.
To this end we will need to determine concurrently the sequence: 
\eq{eq:DGLpexp} 
\psi_k = \ \int_{-1}^{1} L_k'(s) e^{- \half ((s-t)/\sigma) ^2} ds .
\en
From the 3-term recurrence of the Legendre polynomials:
\eq{eq:Leg}
(k+1) L_{k+1} (\lambda) = (2 k+1) \lambda L_k(\lambda) - k L_{k-1} (\lambda) 
\en
we get by integration:
\eq{eq:Leg1} 
(k+1) \gamma_{k+1} = (2 k+1) \int_{-1}^{1} s L_k(s) 
e^{- \half ((s-t)/\sigma) ^2} ds 
- k \gamma_{k-1}   .
\en 
A useful observation is that the above formula is valid for $k=0$ by setting
$\gamma_{-1} \equiv 0$. This comes from  \nref{eq:Leg}, which
is valid for $k=0$ by setting $L_{-1} (\lambda) \equiv 0$. 
Next we expand the integral term in the above equality:
\begin{eqnarray} 
\ \int_{-1}^{1} s  e^{- \half ((s-t)/\sigma) ^2}  L_k(s) ds 
& = & \sigma^2 \ \int_{-1}^{1} \frac{ s - t}{\sigma^2}   
e^{- \half ((s-t)/\sigma) ^2} L_k(s) ds 
+ t  \gamma_k \\ 
& = & \sigma^2 \ \int_{-1}^{1} \frac{ d }{ds} [-  e^{- \half ((s-t)/\sigma)^2}]
 L_k(s) ds + t  \gamma_k.
\label{eq:gam22} 
\end{eqnarray} 
The next step is to proceed with integration by parts for the integral in
the above expression: 
\begin{eqnarray} 
\int_{-1}^{1} \frac{ d }{ds} [-  e^{- \half ((s-t)/\sigma)^2}]
 L_k(s) ds &=&  -L_k (s) e^{- \half ((s-t)/\sigma)^2} \Big\vert_{-1}^1  \nonumber \\ 
&+&
\int_{-1}^{1} e^{- \half ((s-t)/\sigma)^2} L_k'(s) ds .
\label{eq:gam23} 
\end{eqnarray} 
Noting that $L_k(1) = 1$ and $L_k(-1) = (-1)^k$ for all $k$, we get
\begin{small}
\begin{eqnarray} 
\int_{-1}^{1} \frac{ d }{ds} [-  e^{- \half ((s-t)/\sigma)^2}]
 L_k(s) ds 
& =   &
- e^{- \half ((1-t)/\sigma)^2} 
+ (-1)^k  e^{- \half ((1+t)/\sigma)^2} 
+ \psi_k  
\label{eq:gam24} \\ 
& =   & -  e^{-\half (1+t^2)/\sigma^2}  
\left[ e^{t/\sigma^2} - (-1)^k  e^{-t/\sigma^2}\right] + \psi_k 
\label{eq:gam245} \\
& \equiv & \psi_k - \zeta_{k}  , \label{eq:gam25} 
\end{eqnarray} 
\end{small}
where we have defined 
\eq{eq:eta} 
\begin{split}
\psi_{k} &=
\int_{-1}^{1} e^{- \half ((s-t)/\sigma)^2} L_k'(s) ds,\\
\zeta_{k} &= 
 e^{- \half ((1-t)/\sigma)^2} 
- (-1)^k  e^{- \half ((1+t)/\sigma)^2}.
\end{split}
\en
We note in passing that according to \nref{eq:gam245}, 
$\zeta_k$ can be written as
\[
\zeta_{k} = \left\{ \begin{array}{ll} 
2 e^{-\half (1+t^2)/\sigma^2}  \mbox{sh} 
(t/\sigma^2 ) & \mbox{for $k$ even} \\
2 e^{-\half (1+t^2)/\sigma^2}  \mbox{ch} 
(t/\sigma^2 ) & \mbox{for $k$ odd} . 
\end{array} \right.
\]
Substituting \nref{eq:gam24} into \nref{eq:gam22} and the result into 
\nref{eq:Leg1} yields
\eq{eq:Leg2}  
(k+1) \gamma_{k+1} = 
(2k+1) \left[\sigma^2(\psi_k - \zeta_{k}) + t \gamma_k \right] -k \gamma_{k-1}
\en 

The only thing that is left to do is to find a recurrence for the 
$\psi_k$'s. Here we use the elegant formula which can be found in, e.g., 
\cite[p. 47]{Lebedev-book}
\eq{eq:Lp}
L_{k+1}' (\lambda) = (2k+1) L_k (\lambda) + L_{k-1}' (\lambda) .
\en
Integrating in $[-1, 1 ]$  yields the relation: 
\eq{eq:psik1}
\psi_{k+1} = (2k+1) \gamma_k + \psi_{k-1} 
\en
Note that initial values of $\psi_k $ are
$\psi_0 = 0$, $\psi_1  = \gamma_0$.
In the end, we obtain the following recurrence relations: 
\eq{eq:Leg3} 
\left\{ \begin{array}{lll} 
\gamma_{k+1}  & = &
\frac{2k+1}{k+1}
 \left[\sigma^2(\psi_k-\zeta_{k}) + t \gamma_k \right] - 
\frac{k}{k+1}  \gamma_{k-1}\\
\psi_{k+1} & = & (2k+1) \gamma_k + \psi_{k-1} . 
\end{array} 
\right. 
\en 
It can be noted that the above formulas work for $k=0$ by setting
$\gamma_{-1} = \psi_{-1} = 0$. 
The recurrence starts with $k=0$, using the initial values
$\gamma_0$ given by \nref{eq:gam0}, 
$\psi_1 = \gamma_0$, and $\psi_0 = 0$. 

An important remark here is that one has to be careful about
the application of the recurrence \nref{eq:Leg3}. The perceptive
reader may notice that such a recurrence runs the risk of being
unstable. In fact we observe the following behavior. For 
large values of $\sigma$ the Gaussian function can be very smooth and 
as a result a very small degree of polynomials may be needed, i.e., the
value of $\gamma_k$ drop to small values quite rapidly as $k$ increases.
If we ask for a high degree polynomial 
and continue the recurrence \nref{eq:Leg3}  beyond the point where the
expansion has converged (indicated by small value of $\gamma_k$)
we will essentially iterate with noise. As it turns out, this noise is 
amplified by the recurrence. This is because the coefficient 
$\psi_k - \zeta_k $ becomes just noise and this causes the recurrence to
diverge. An easy remedy is to just stop iterating \nref{eq:Leg3} as soon 
as two consecutive $\gamma_k$'s are small. This takes care of two issues
at the same time. First, it determines a sort of optimal degree to be
used. Second, it avoids the unstable behavior observed by continuing the
recurrence. Specifically, a test such as the following is performed:
\eq{eq:StopRecur}
|\gamma_{k-1}| + |\gamma_{k}| \le k \cdot \text{tol},
\en
where $\text{tol}$ is a small tolerance, which can be 
set to $10^{-6}$ for example.

With this we can now easily develop the Delta-Gauss-Legendre (DGL) expansion algorithm.
In the DGL algorithm, we will refer to formula \nref{eq:DelDens}. But now $p_M$ 
is the $M$-degree polynomial 
\eq{eq:pmnew} 
p_M(\lambda) = \frac{1}{(2\pi \sigma^2)^{1/2}} 
\sum_{k=0}^M \left( k + \frac{1}{2} \right)  \gamma_k  L_k(\lambda) ,
\en 
obtained by truncating the sum \nref{eq:DGLexp} to $M+1$ terms. 

\section{Cumulative density of states from the Lanczos
method}\label{sec:app_lanczos}

\begin{figure}[ht]
	\begin{center}
		\includegraphics[width=0.35\textwidth,height=0.35\textwidth]{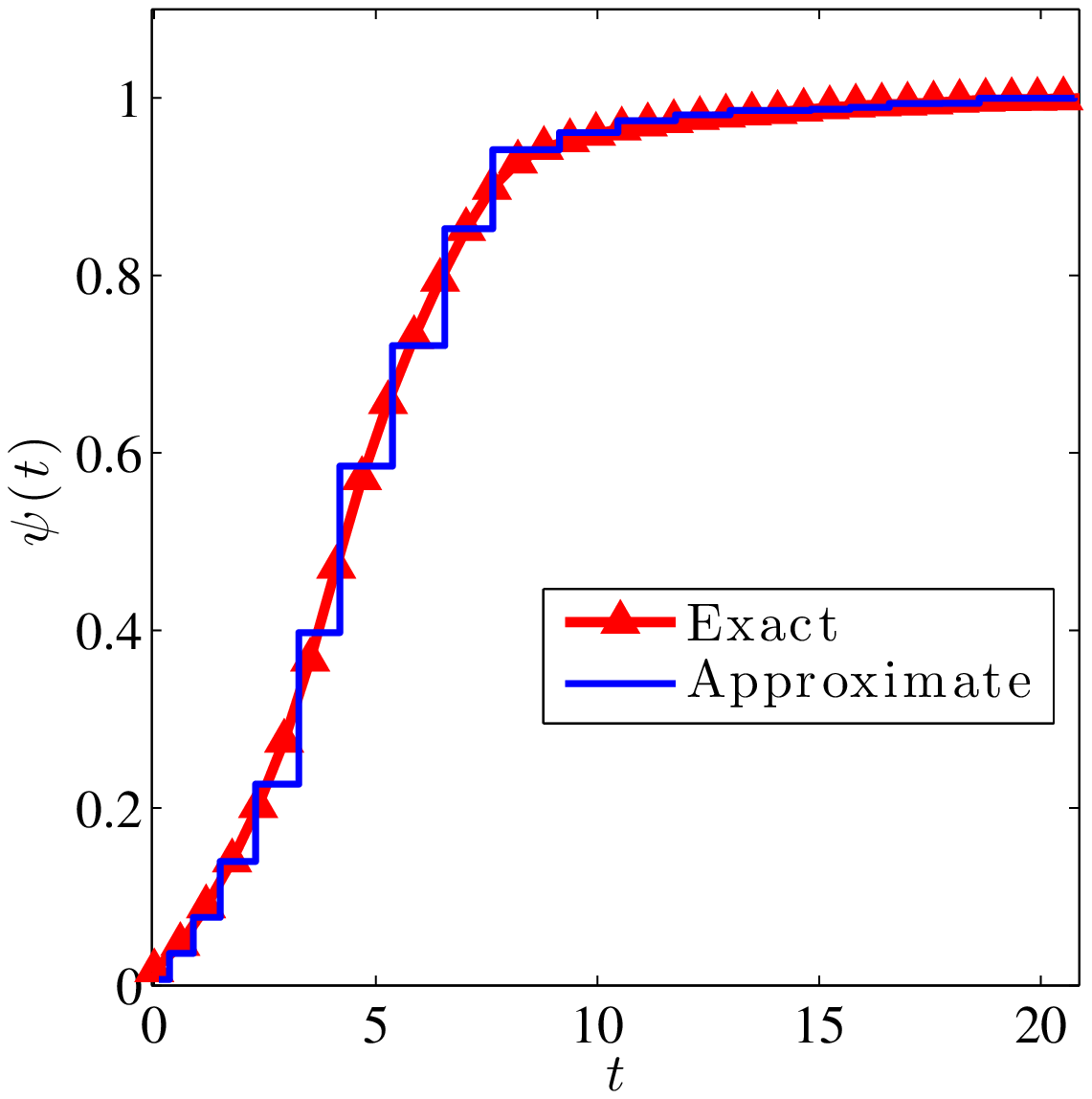}\quad
		\includegraphics[width=0.35\textwidth,height=0.35\textwidth]{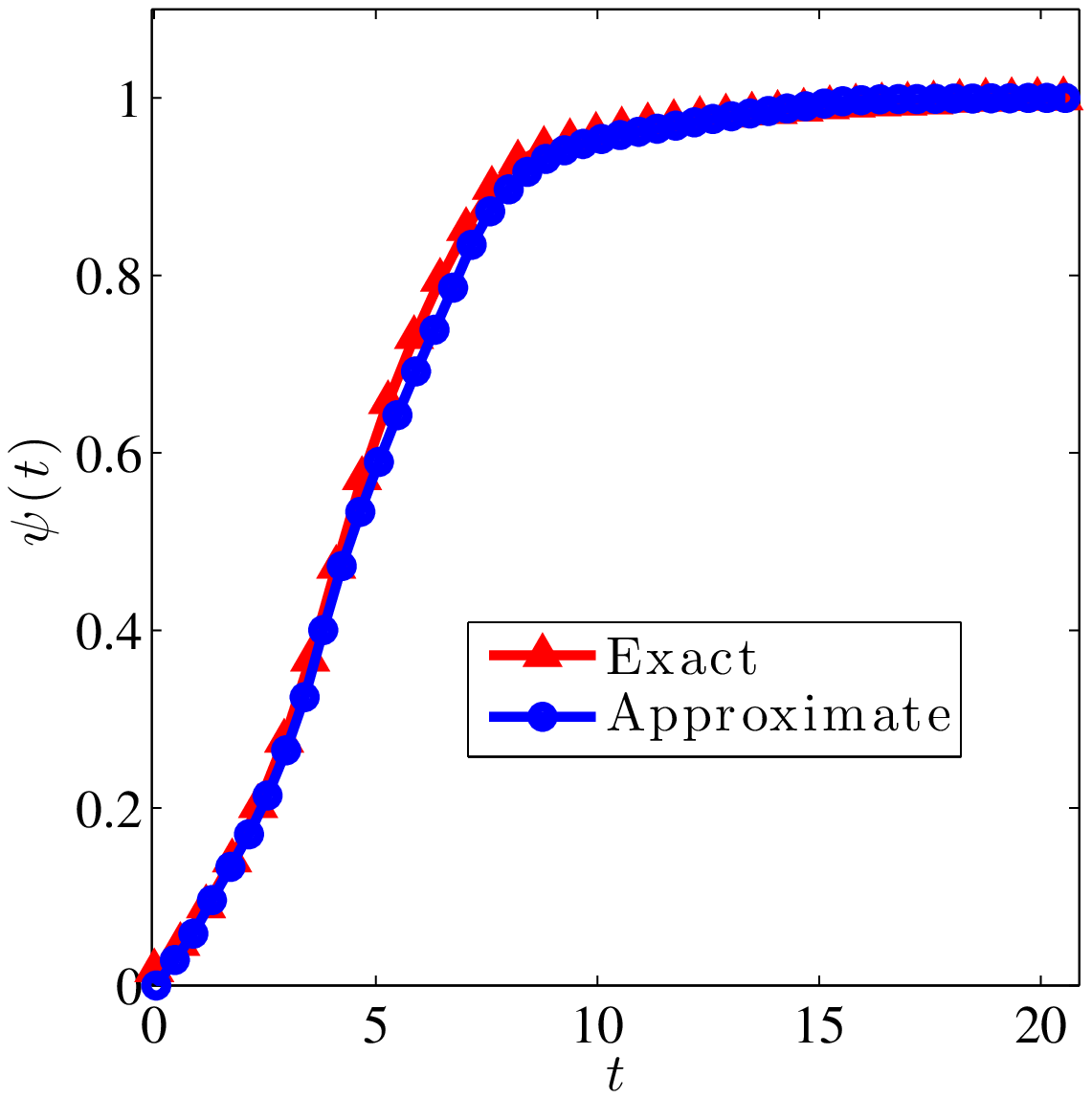}
	\end{center}
	\caption{The approximate cumulative spectral density associated with
	the modified 2D Laplacian constructed directly from a 20-step Lanczos
	run (left) and its spline-interpolated and smooth version (right).}
	\label{fig:lancdos}
\end{figure}

An alternative way to refine 
the Lanczos based DOS approximation from a $M$-step
Lanczos run is to first construct an approximate cumulative spectral density 
or cumulative density of states (CDOS), defined as
\[
\psi(t) = \int_{\infty}^{t}\phi(s) ds.
\]
Without applying regularization,
the approximate CDOS can be computed from the Lanczos procedure as
\begin{equation}
\wt{\psi}(t) =  
\sum_{k = 0}^M \eta_{k}^2 \delta(t-\theta_{k}),
\label{lancdos}
\end{equation}
where $\eta_{k}^2 = \sum_{i=1}^k \tau_{i}^2$, and $\theta_{k}$ and 
$\tau_{k}$ are eigenvalues and the first components of the 
eigenvectors of the tridiagonal matrix $T_{M}$ defined in~\eqref{eq:lanapprox}. 
This approximation is plotted as a staircase function in 
Figure~\ref{fig:lancdos} for the modified 2D Laplacian.
%
%
Note that both $\psi(t)$ and $\wt{\psi}(t)$ are monotonically non-decreasing functions.
Furthermore, it can been shown \cite{ks53,FischerFreundW1994} that $\psi(t)-\wt{\psi}(t)$ has 
precisely $2M-1$ sign changes within the spectrum of $A$. A sign change
occurs when $\psi(t)$ crosses either a vertical or horizontal step 
of $\wt{\psi}(t)$.  These properties allow us to construct an 
``interpolated'' CDOS that matches $\psi(t)$ and $\wt{\psi}(t)$ at the 
points where $\psi(t)$ crosses $\wt{\psi}(t)$.

\bibliographystyle{siam}
\bibliography{paper}


\end{document}

%% file: preamble.tex
\def\myrefs#1#2{ 
{\bigskip \noindent
{\Large \bf #2}  
 \list {[\arabic{enumi}]}{\settowidth\labelwidth{[#1]}
 \leftmargin\labelwidth 
 \advance\leftmargin\labelsep
 \usecounter{enumi} }  
 \def\newblock{\hskip .11em plus .33em minus .07em}
 \sloppy\clubpenalty4000\widowpenalty4000
 \sfcode`\.=1000\relax}  }

\def\half{{1\over2}}%
\def\nref#1{(\ref{#1})}

\def\comb#1,#2,{ \left( {#1 \atop #2 } \right)  }%
\def\prodd#1,#2,#3,{ \prod_{\scriptstyle #1 \atop\scriptstyle #2 }^{#3} }%
\def\summ#1,#2,#3,{ \sum_{\scriptstyle #1 \atop\scriptstyle #2 }^{#3} }%
\newtheorem{algor}{{\sc Algorithm}}[section]
\newtheorem{tabl}{Table}[section]

\def\betab{\begin{tabbing} 
xxxx\=xxxx\=xxx\=xx\=xx\=xx\=xx\=xx\=xx\=xx\=xx\=xx\=xx\= \kill} 
\def\entab{\end{tabbing}\vspace{-0.12in}}

\newcommand{\eq}[1]{\begin{equation}\label{#1}}
\newcommand{\en}{\end{equation}}

\newcommand{\beeq}[1]{\begin{equation}\label{#1}}
\newcommand{\eneq}{\end{equation}}

\renewcommand{\Im}{\mathrm{Im}~}

\newcommand{\abs}[1]{\lvert#1\rvert}

\newcommand{\average}[1]{\langle#1\rangle}

\newcommand{\nvec}{n_{\mathrm{vec}}}
\newcommand{\wt}{\widetilde}